\documentclass{article}

\usepackage{amssymb}
\usepackage{amsmath}
\usepackage{amsthm}

\usepackage{amsfonts}
\usepackage{xcolor}
\usepackage{fullpage}
\usepackage{hyperref}
\usepackage{algorithm}
\usepackage{algpseudocode}
\usepackage{graphicx, caption, subcaption}
\usepackage{framed}
\usepackage{booktabs}
\usepackage{float}

\DeclareMathOperator{\rank}{rank}

\newtheorem{theorem}{\bf Theorem}[section]

\newtheorem{lemma}[theorem]{\bf Lemma}

\newtheorem{assumption}{Assumption}

\newcommand{\R}{\mathbb{R}}

\newcommand{\Rnn}{\mathbb{R}^{n \times n}}
\newcommand{\Rmn}{\mathbb{R}^{m \times n}}

\newcommand{\Rnr}{\mathbb{R}^{n \times r}}

\newcommand{\Rmr}{\mathbb{R}^{m \times r}}

\newcommand{\Rrr}{\mathbb{R}^{r \times r}}

\newcommand{\norm}[1]{\left\| #1 \right\|}

\newcommand{\Mr}{\mathcal{M}_r}

\newcommand{\Tr}{\mathcal{T}_r}

\newcommand{\proj}[2]{\mathcal{P}_{#1} \left[ #2 \right]}

\newcommand{\Range}{\mathrm{Range}}

\usepackage{accents}
\newcommand*{\dt}[1]{%
  \accentset{\mbox{\large\bfseries .}}{#1}}

\title{Randomized methods for dynamical low-rank approximation}

\author{Benjamin Carrel}

\begin{document}

\maketitle

\begin{abstract}
We introduce novel dynamical low-rank methods for solving large-scale matrix differential equations, motivated by algorithms from randomized numerical linear algebra.
In terms of performance (ratio accuracy/cost), our methods can overperform existing dynamical low-rank techniques.
Several applications to stiff differential equations demonstrate the robustness, accuracy and low variance of the new methods, despite their inherent randomness. Allowing augmentation of the range and corange, the new methods have a good potential for preserving critical physical quantities such as the energy, mass and momentum. Numerical experiments on the Vlasov-Poisson equation are particularly encouraging. 

The new methods comprise two essential steps: a range estimation step followed by a post-processing step.
The range estimation is achieved through a novel dynamical rangefinder method.
Subsequently, we propose two methods for post-processing, leading to two time-stepping methods: dynamical randomized singular value decomposition (DRSVD) and dynamical generalized Nyström (DGN).
The new methods naturally extend to the rank-adaptive framework by estimating the error via Gaussian sampling.
\end{abstract}

\section{Introduction} \label{sec: intro}

In this paper, we are interested in solving autonomous matrix differential equations of the general form
\begin{equation} \label{eq: full order model}
\dt{A}(t) = F(A(t)), \quad A(0) = A_0 \in \Rmn,
\end{equation}
for some vector field $F: \Rmn \rightarrow \Rmn$ and time interval $[0,T]$.
Examples of such differential equations include the discretization of partial differential equations in physics~\cite{einkemmer2018low,cassini2022efficient,coughlin2022efficient,einkemmer2024accelerating}, uncertainty quantification~\cite{babaee2017robust,ali2024dynamically,feppon2018dynamically}, quantum chemistry~\cite{sulz2024numerical,ceruti2023rank} and also machine learning~\cite{schotthofer2022low,savostianova2023robust,schmidt2023rank}; further applications will be provided in the next paragraphs.
We refer to equation~\eqref{eq: full order model} as the full order model, and when the parameters $m$, $n$ are large, computing a dense solution to the problem becomes both computationally and memory intensive and is sometimes not feasible with current hardware.
To address this issue, reduced order modeling techniques have grown in popularity. 
A particularly successful one is the dynamical low-rank approximation~\cite{koch2007dynamical,einkemmer2018low}, which uses low-rank compressions to efficiently store and evolve the solution over time.

On the other hand, the domain of randomized numerical linear algebra has emerged and is growing in popularity. A key achievement in the domain is the rangefinder algorithm, also called HMT, which comes with good theoretical guarantees~\cite{halko2011finding}. A direct application of the rangefinder is the randomized singular value decomposition (randomized SVD) algorithm which allows to approximate the truncated SVD of a (large) matrix $A \in \Rmn$ at a reduced cost. The randomized SVD algorithm scales well since it only requires a few matrix-vector multiplications, and such operations are particularly fast on GPUs.
The generalized Nyström (GN) scheme~\cite{woodruff2014sketching,tropp2017practical,nakatsukasa2020fast} is also a fast method for approximating the SVD of low-rank matrices, it is numerically stable and also more accurate than the randomized SVD.

In this paper, we propose new time-stepping methods to approximate equation~\eqref{eq: full order model} by means of low-rank approximations. The new methods' construction is inspired from the rangefinder, the randomized SVD and the generalized Nyström. At the moment, we are not aware of a direct formulation of randomized methods in the context time-dependent problems of the form~\eqref{eq: full order model}, this paper is therefore an original attempt to link the domain of dynamical low-rank approximation and the domain of randomized numerical linear algebra. The new methods show a good potential for various applications in physics, uncertainty quantification, quantum chemistry, and machine learning. Before defining the new methods, we introduce notations and key results in the two domains.

\subsection{Preliminary on randomized numerical linear algebra}

\subsubsection*{The rangefinder algorithm}

The idea of using a randomized subspace iteration for finding the range goes back to \cite{rokhlin2010randomized}. 
Quickly after, the \textit{powered randomized rangefinder} algorithm was proposed and analyzed in \cite{halko2011finding}, and the algorithm is now simply called rangefinder or HMT. This algorithm solves the problem of approximating the range of a matrix $A \in \Rmn$ with a prescribed target rank $r$ and oversampling parameter $p$. More precisely, the rangefinder computes a matrix $Q \in \R^{m \times (r+p)}$ with orthonormal columns such that
$$\norm{A - Q Q^T A}_2 = \min!$$
The rangefinder algorithm consists of the following steps:
\begin{enumerate}
  \item Draw a random Gaussian matrix $\Omega \in \R^{n \times (r+p)}$.
  \item Compute the product: $B = A \Omega$.
  \item Extract the range: $Q = \mathrm{orth}(B)$.
\end{enumerate}
The algorithm scales particularly well since it only requires the application of $A$ to $r+p$ random vectors. The range is then extracted from the resulting tall and skinny $m \times (r+p)$ matrix $B$.
The approximation can be improved by performing a few power iterations, which we summarize as follows:
\begin{equation} \label{eq: rangefinder maths}
Q = \mathrm{orth}\left( (AA^T)^q A \Omega \right),
\end{equation}
and numerical stability is achieved with an QR step after each application of $A$.
The detailed procedure is given in~\cite[Algorithm 4.4]{halko2011finding}, and we recall here the main convergence theorem from~\cite[Corrolary 10.10]{halko2011finding}.
\begin{theorem}[Average spectral error of the rangefinder] \label{thm: rangefinder}
  Let $A \in \Rmn$ with singular values $\sigma_1 \geq \sigma_2 \geq \ldots$. 
  Choose a target rank $r \geq 2$ and an oversampling parameter $p \geq 2$, where $r+p \leq \min \{m,n\}$. 
  Draw an $n \times (r+p)$ random Gaussian matrix $\Omega$ and construct $Q = \mathrm{orth}\left( (AA^T)^q A \Omega \right)$.
  Then, the expected approximation error verifies
  $$\mathbb{E} \norm{A - Q Q^T A}_2 \leq \left[ \left( 1 + \sqrt{\frac{r}{p-1}} \right) \sigma_{r+1}^{2q+1} + \frac{e \sqrt{r + p}}{p} \left( \sum_{j > r} \sigma_j^{2(2q+1)} \right)^{1/2} \right]^{1/(2q+1)}.$$
\end{theorem}

In practice, the rangefinder is typically used with a relatively small oversampling parameter, typically $p~\leq~r$. With power iterations, the constant in front of the best approximation error $\sigma_{r+1}$ contracts exponentially fast to one. While random Gaussian matrices are generally a good choice, it appears that sub-sampled randomized trigonometric transforms (SRTTs) are slightly more accurate~\cite{martinsson2020randomized}. Alternatively, sparse maps allow faster matrix-vector multiplications.

\newpage

\subsubsection*{The randomized SVD algorithm}

A direct application of the rangefinder algorithm is the randomized SVD algorithm \cite{halko2011finding} which we briefly summarize below.
\begin{enumerate}
  \item Estimate the range by computing $Q = \mathrm{Rangefinder}(A, r, p, q)$.
  \item Compute the reduced SVD of an $(r+p) \times n$ matrix: $\tilde U \Sigma V^T = \mathrm{svd}(Q^T A)$.
  \item Update the basis and truncate the approximation: $\Tr(A) \approx \Tr(Q \tilde U \Sigma V^T)$.
\end{enumerate}
The operator $\mathcal T_r : \Rmn \rightarrow \Rmn$ is the SVD truncated to rank $r$, also denoted $\mathcal T_r(A) = A_{[r]}$ in the literature.
The randomized SVD scales particularly well since it only requires a few matrix-vector multiplications, and then costly operations are performed on \textit{small} rectangular matrices. 
The truncation step is cost-free since the input is already given by its SVD, which can be truncated by keeping the $r$ first singular vectors/values.
The analysis follows directly from Theorem~\ref{thm: rangefinder}, see also~\cite[Theorem 10.5]{halko2011finding}.

Since then, the method became standard for computing the SVD of large-scale matrices, and has also been extended to a variety of problems, including the principal component analysis \cite{rokhlin2010randomized,greenacre2022principal}, matrix completion \cite{feng2018faster,wen2023accelerated}, kernel matrix approximation \cite{li2014large}, generalized SVD \cite{saibaba2021randomized}, dynamic mode decomposition \cite{schmid2022dynamic} regularization of inverse problems \cite{xiang2013regularization}, image compression with GPU \cite{ji2014gpu}, tensor structures \cite{batselier2018computing,ahmadi2021randomized}, and many more.

\subsubsection*{Generalized Nyström schemes}

Generalized Nyström (GN) schemes are an alternative to the randomized SVD algorithm.
Let us introduce first Nyström's scheme, originally proposed in \cite{nystrom1930praktische} and more recently popularized in \cite{williams2000using}.
Let $A~\in~\Rnn$ be a symmetric square positive semi-definite (SPSD) matrix. 
Then, a rank-$r$ approximation is given by
$$A \approx AX (X^T A X)^{\dagger} (AX)^T,$$ 
where $X$ is an $n \times r$ random sketch matrix. 
The scheme became popular in machine learning for working with kernel matrices \cite{williams2000using,scholkopf2002learning,fowlkes2004spectral}, and has been extensively analyzed in \cite{gittens2011spectral,gittens2013revisiting}.
Interestingly, it has been observed in \cite{halko2011finding,martinsson2019randomized} and studied in \cite{gittens2013revisiting} that computing
$$A \approx AQ (Q^T A Q)^{\dagger} (AQ)^T, \quad Q = \mathrm{Rangefinder}(AX)$$
outperforms the randomized SVD in most situations.
Later on, Nyström's scheme has been generalized to rectangular matrices.
Early algorithms were proposed in \cite{clarkson2009numerical,woodruff2014sketching,tropp2017practical}, and we refer to \cite{nakatsukasa2020fast} for a complete description and analysis of the scheme.
For the purpose of the paper, we are interested in a particular variant of the \textit{generalized Nyström} algorithm which uses the rangefinder algorithm and a second oversampling parameter $\ell \geq0$, and which we define as follows:
\begin{equation} \label{eq: orthogonal generalized Nyström}
  A \approx \tilde A = (A W) \ \Tr(Q^T A W)^{\dagger} \ (Q^T A), \qquad Q = \mathrm{Rangefinder}(A, r, p, q), \quad W = \mathrm{Rangefinder}(A^T, r, p+\ell, q).
\end{equation}
As pointed out in \cite{nakatsukasa2020fast}, the order in which operations are done is important for the numerical stability of the algorithm. Following our notation, the truncation to the target rank $r$ is applied first, and then the pseudoinverse of the truncated matrix is computed. Computing the pseudoinverse of an SVD is cheap since $(U \Sigma V^T)^{\dagger}~=V\Sigma^{-1}U^T$.
Overall, the cost of this version of generalized Nyström is roughly double the cost of a randomized SVD.

\subsection{Dynamical low-rank approximation}

The dynamical low-rank approximation (DLRA) was originally proposed in \cite{koch2007dynamical} and in \cite{nonnenmacher2008dynamical} as a technique for evolving problems using low-rank factorizations to compress data. 
Let us now introduce important notations, briefly recall the fundamentals and then do a quick overview of the literature on the topic. In a separated paragraph below, we further discuss the differences between the new methods proposed in this paper and existing DLRA integrators.

\subsubsection*{Notations and theory}

We consider the set of fixed rank matrices
$$\mathcal M_r = \{M \in \R^{m \times n} \mid \mathrm{rank}(M)=r \},$$
which is in fact a smooth manifold~\cite{lee2003smooth} embedded into $\Rmn$. 
The dynamical low-rank approximation is essentially the time-dependent variational principle (TDVP) applied to the vector field of the full order model. Restricted to the above manifold of low-rank matrices, the DLRA evolves the projected dynamics
\begin{equation} \label{eq: dlra model}
\dt{Y}(t) = \mathcal P_{Y(t)} F(Y(t)), \quad Y(0) = Y_0 \in \Mr,
\end{equation}
where $\mathcal P_{Y}$ denotes the $\ell_2$-orthogonal projection onto the tangent space $\mathcal T_Y \Mr$, here defined at a certain low-rank matrix $Y~\in~\Mr$.
The analysis of the accuracy of the DLRA relies on the following two assumptions.
\begin{assumption}[One-sided Lipschitz] \label{ass: one sided Lipschitz}
  There exists a constant $\ell \in \R$ such that for all $X, Y \in \Rmn$ we have
  $$\langle X - Y,\, F(X) - F(Y) \rangle \leq \ell \norm{X - Y}_F^2.$$
\end{assumption}
\begin{assumption}[Closeness to the tangent bundle] \label{ass: closeness}
  There exists a small constant $\varepsilon_r > 0$ and a suitable neighborhood~$\mathcal N$ of the solution to the full order model~\eqref{eq: full order model} such that for all $Y \in \Mr \cap \mathcal N$ we have
  $$\norm{F(Y) - \proj{Y}{F(Y)}}_F \leq \varepsilon_r.$$
\end{assumption}
Note that the one-sided Lipschitz constant can be negative, e.g. on diffusive problems.
Under these two assumptions, the solution $Y(t)$ to the DLRA compared to the full order model solution $A(t)$ verify the error bound
\begin{equation} \label{eq: dlra error}
\norm{A(t) - Y(t)}_F \leq e^{t \ell} \norm{A_0 - Y_0}_F + \varepsilon_r \int_0^t e^{(t-s) \ell} ds,
\end{equation}
and we refer to \cite{kieri2019projection} for a proof.
The DLRA model given in \eqref{eq: dlra model} is equivalent to a system of three coupled differential equations~\cite{koch2007dynamical}. Solving the coupled system is numerically ill-conditioned in the presence of small singular values, and therefore robust techniques have been proposed in the literature.

\subsubsection*{A quick overview of numerical integrators}

As above mentioned, solving the projected differential equation~\eqref{eq: dlra model} can be numerically ill-conditioned in the presence of small singular values. Remember that we do expect to have small singular values since we are interested in problems with good low-rank approximations. It is therefore crucial to have methods that are numerically robust to the stiffness induced by small singular values. We now give a short overview of the literature on robust numerical methods and some of their applications. 

The projector-splitting integrator~\cite{lubich2014projector} was the first integrator proposed, and its robust convergence was proved in~\cite{kieri2016discretized}. 
To overcome the limitations of the backward step in the projector-splitting integrator, bases updates and Galerkin (BUG) methods were proposed in~\cite{ceruti2022unconventional,ceruti2022rank,ceruti2024parallel,ceruti2024robust,kusch2024second}. 
Higher-order methods based on low-rank explicit Runge--Kutta schemes were proposed in~\cite{kieri2019projection,lam2024randomized,nobile2025robust}, and low-rank exponential Runge--Kutta methods for stiff differential equations are studied in~\cite{carrel2023projected}. 
Low-rank implicit schemes have also been studied in~\cite{nakao2025reduced,yin2023semi,appelo2025robust,sutti2024implicit,rodgers2023implicit}, and a parallel-in-time low-rank integrator based on the parareal algorithm is discussed in~\cite{carrel2023low}. Dynamical CUR decomposition techniques have also been proposed in~\cite{donello2023oblique,naderi2024cur,park2025low}, these methods have often less theoretical guarantees but have shown excellent results in practice and seem to be easier to implement and computationally less intensive than DLRA techniques, even though being closely related.
Many of these methods have natural extensions to higher-dimensional problems with imposed tensor structures such as Tucker, tensor trains (or matrix product states) and tree tensor networks, see for example~\cite{koch2010dynamical,lubich2013dynamical,lubich2018time,ceruti2020time,ceruti2022unconventional,ceruti2023rank,ghahremani2024cross,einkemmer2024hierarchical,ceruti2024parallelTTNs}.

The literature on applications of the dynamical low-rank approximation is now quite extensive, so we try to give a quick overview and more references can be found within the following papers. The DLRA has been successfully applied to various problems arising in quantum chemistry~\cite{sulz2024numerical,ceruti2024parallelTTNs} and in physics: kinetics equations~\cite{einkemmer2018low,einkemmer2019low,einkemmer2020low}, radiation therapy~\cite{kusch2023robust}, the Su--Olson problem~\cite{baumann2024energy}, fission criticality~\cite{scalone2025multi} and also to stochastic differential equations~\cite{kazashi2024dynamical,nobile2025petrov,nobile2023error}. Key properties of the numerical integrators in the context of physics are studied in~\cite{einkemmer2021mass,einkemmer2023robust,kusch2023stability,einkemmer2023conservation}. Recent research also show applications in machine learning, such techniques are called dynamical low-rank training (DLRT) and are studied in~\cite{schotthofer2022low,savostianova2023robust,schotthofer2024federated,schotthofer2024geolora,schotthofer2025dynamical,kusch2025augmented,billaud2017dynamical,billaud2022new}.

\subsection{Contributions and outline} \label{sec: contributions}

As discussed in the previous section, the numerical issues due to the small singular values have been addressed by several numerical integrators. However, the low-rank integration of stiff differential equations is still challenging as they impose strong stepsize restrictions. Robust integrators based on exponential methods have been proposed in~\cite{ostermann2019convergence,carrel2023projected}, but they are limited to stiff differential equations with a Sylvester-like structure, and little is known for the low-rank integration of stiff differential equations with more general structures. The analysis of BUG-based methods \cite{ceruti2022unconventional,ceruti2024parallel,kusch2024second} and projected Runge--Kutta methods \cite{kieri2019projection} rely on the Lipschitz constant of the function of the vector field, implying large hidden constants which can be numerically observed.

The idea of bringing randomization to the dynamical low-rank approximation has been recently studied in~\cite{lam2024randomized} where they propose randomized low-rank Runge--Kutta methods. The key idea of these methods is to use randomization (with the generalized Nyström formula) within explicit Runge--Kutta schemes to limit the rank growth in the stages. By doing so, they avoid the computation of projections onto the tangent space, which makes the schemes robust to the presence of small singular values, and in some cases more robust than the projected Runge--Kutta methods proposed in~\cite{kieri2019projection}. However, the approach is currently limited to explicit Runge--Kutta schemes and, as a consequence, the convergence depends again on the Lipschitz constant of the vector field, leading to large hidden constants when applied to stiff differential equations.
In contrast to the projector-splitting and BUG integrators, low-rank explicit Runge--Kutta methods do not present a decomposition into substeps, which reduces the potential for parallelization.

Our main motivation is to develop BUG-like methods that can be applied to stiff differential equations without imposing strong stepsize restrictions. To achieve this, we introduce two new randomized time-stepping methods with a structure similar to existing BUG-based integrators. As a consequence, the new methods have a good potential for parallelization (in particular the dynamical generalized Nyström) and existing BUG implementations can be easily modified to implement the new methods. Additionally, we present novel rank-adaptive time-stepping methods where a Gaussian sampling is performed to estimate the errors at little extra cost.
Our numerical experiments indicate that the new methods often outperform existing DLR integrators in terms of the ratio between accuracy and cost, and demonstrate particular effectiveness in handling stiff differential equations characterized by large Lipschitz constants. With augmentation of the range and corange, the new methods also have a good potential for preserving critical physical quantities. A comprehensive convergence analysis of these methods remains an avenue for future research.

The paper is organized as follows. In Section~\ref{sec: dynamical rangefinder}, we propose a new method that we call dynamical rangefinder, which can be used to estimate the range (or corange) of the solution of a matrix differential equation with a prescribed rank, and we also propose an alternative rank-adaptive method using a prescribed tolerance instead. We provide numerical experiments on a toy problem to show the accuracy of the estimated range, and we compare the new method with the rangefinder directly applied to the reference solution.
Next, we build on top of the new dynamical rangefinder, and we propose in Section~\ref{sec: time stepping mehtods} two new time-stepping methods that we respectively call dynamical randomized SVD and dynamical generalized Nyström.
We show that, when the solution to the full order model is exactly low-rank and that the estimated range contains the range of the full solution, then the new methods are exact.
In Section~\ref{sec: applications}, we apply the two new time-stepping methods to various problems on which we compare the performance (ratio accuracy/cost) with different kind of DLRA techniques.

\section{Dynamical randomized rangefinder} \label{sec: dynamical rangefinder}

Recall the full order problem \eqref{eq: full order model} defined by the differential equation
$$\dt{A}(t) = F(A(t)), \quad A(0) = A_0 \in \Rmn.$$ 
A common notation for the numerical flow computed by a numerical solver is $\phi^F_t(A_0)$ and $\varphi^F_t(A_0)$ for the exact flow (or exact solution). In this paper, we are not studying the error of numerical flows but only considering exact flows. Since the dependence on the initial value is not relevant for us, we use the explicit notation $A(t)$ for denoting the exact flow of the full order problem at time $t$, that is the exact solution to the full order model.

Let us fix a stepsize $h \in (0,T]$ and a target rank of approximation $r$, we now focus on the problem of finding the range which approximates best the exact solution $A(h)$. In other words, we seek a matrix $Q_h \in \R^{m \times r}$ such that
$$\norm{A(h) - Q_h Q_h^T A(h)}_2 = \min !$$
Here, the notation "$\min!$" is similar to the notation used in \cite{feischl2024regularized}.
In this section, we start by discussing an idealized method where the exact solution to the full order model is available. Then, we propose a practical method to find such a matrix $Q_h$ without having to solve the full order model, and we propose a way to improve the approximation with dynamical power iterations. Finally, we discuss a method for automatically choosing the approximation rank for a prescribed tolerance $\tau$ and probability of failure $\delta$.

\subsection{An idealized dynamical rangefinder method} \label{sec: idealized dynamical rangefinder}

As it is sometimes done in the literature on the dynamical low-rank approximation (e.g. in \cite{lubich2014projector} and \cite{ceruti2022unconventional}), we start by studying the problem on the exact flow $A(t)$, and we consider the following steps:
\begin{itemize}
  \item Sketch a random Gaussian matrix $\Omega \in \R^{n \times (r+p)}$.
  \item Solve the differential equation
  \begin{equation} \label{eq: exact B-step} 
    \left\{
  \begin{aligned}
  &\dt{\mathcal B}(t) = F(A(t)) \Omega, \quad t \in [0,h] \\
  &\mathcal B(0) = A_0 \Omega. 
  \end{aligned}
    \right.
  \end{equation}
  \item Extract the basis at final time: $\mathcal Q_h = \mathrm{orth}(\mathcal B(h))$.
\end{itemize}
When the rank of the solution is over-estimated, the following lemma shows that the above procedure is exact.
\begin{lemma}[Exact range property] \label{lemma: exact range property}
    Fix a sketch matrix $\Omega \in \R^{n \times (r+p)}$ and assume that $\rank(A(h)) \leq r + p$. 
    Then, the matrix $\mathcal Q_h$ obtained from the above procedure verifies $\Range( \mathcal Q_h) = \Range(A(h))$ and therefore
    $$\norm{A(h) - \mathcal Q_h \mathcal Q_h^T A(h)}_2 = 0.$$
\end{lemma}
\begin{proof}
    Suppose that $\mathrm{rank}(A(h)) = R \leq r+p$, then we can decompose 
    $$A(h) \Omega = U_h \Sigma_h V_h^T \Omega,$$ 
    where $U_h \in \R^{m \times R}$ and $V_h \in \R^{n \times R}$ have orthonormal columns. 
    Since $\Omega \in \R^{n \times (r+p)}$ is a Gaussian matrix, the result of the product $W =V_h^T \Omega \in \R^{R \times (r+p)}$ is also a Gaussian matrix and has therefore full row-rank with probability one. By definition, we get
    $$\mathcal B(h) = A_0 \Omega + \int_0^h F(A(t)) \Omega \ \mathrm{dt} = A_0 \Omega + \left[ A(h) - A(0) \right] \Omega = A(h) \Omega = U_h \Sigma_h W,$$
    and we conclude by noticing that the following spaces are identical:
    $$\Range(\mathcal Q_h) = \Range(\mathcal B(h)) = \Range(U_h) = \Range(A(h)).$$
\end{proof}
As pointed out in the course of the above proof, the idealized procedure is, in fact, the rangefinder algorithm applied to the exact solution evaluated at time $h$, that is $A(h)$.
The above property is similar to the exactness property of \cite{lubich2014projector} and \cite{ceruti2022unconventional} on exact flows.

\subsection{A practical dynamical rangefinder with dynamical power iterations}

In practice, the exact flow is unknown and equation \eqref{eq: exact B-step} cannot be solved. 
We therefore need to approximate it without solving the full order problem.
To do so, we propose the \textit{dynamical rangefinder} algorithm:
\begin{enumerate}
  \item Sketch a random Gaussian matrix $\Omega \in \R^{n \times (r+p)}$ and compute its pseudoinverse $\Omega^{\dagger} = (\Omega^T \Omega)^{-1} \Omega^T$.
  \item (B-step) Solve the $m \times (r+p)$ differential equation
  \begin{equation} \label{eq: B-step} 
    \left\{
  \begin{aligned}
  &\dt{B}(t) = F(B(t) \Omega^{\dagger}) \Omega, \quad t \in [0,h] \\
  &B(0) = A_0 \Omega. 
  \end{aligned}
  \right.
  \end{equation}
  \item Extract the basis at final time: $Q_h = \mathrm{orth}(B(h))$.
\end{enumerate}
The approximation can further be improved by performing a few power iterations by approximating alternatively the range and corange of the matrix $A(h)$. 
One iteration of the \textit{dynamical power iteration} is done as follows:
\begin{enumerate}
\item[4.] Solve the $n \times (r+p)$ differential equation
  \begin{equation} \label{eq: dynamical rangefinder C-step} 
    \left\{
  \begin{aligned}
  &\dt{C}(t) = F(Q_h C(t))^T Q_h, \quad t \in [0,h] \\
  &C(0) = A_0^T Q_h. 
  \end{aligned}
    \right.
  \end{equation}
\item[5.] Extract the basis at final time: $W_h = \mathrm{orth}(C(h))$. 
\item[6.] Solve the $m \times (r+p)$ differential equation
  \begin{equation} \label{eq: second B-step} 
    \left\{
  \begin{aligned}
  &\dt{B}(t) = F(B(t) W_h^T )  W_h , \quad t \in [0,h] \\
  &B(0) = A_0 W_h. 
  \end{aligned}
    \right.
  \end{equation}
\item[7.]  Extract the basis at final time: $Q_h^{(1)} = \mathrm{orth}(B(h))$.
\end{enumerate}
To perform $q$ dynamical power iterations, one has to repeat steps $4$ to $7$ starting from $Q_h^{(1)}$, then $Q_h^{(2)}$, etc. until obtaining $Q_h^{(q)}$.
The dynamical rangefinder algorithm is summarized in Algorithm \ref{algo: dynamical rangefinder}.
\begin{algorithm}
\caption{Dynamical rangefinder}
\label{algo: dynamical rangefinder}
\begin{algorithmic}[1]
\Require Vector field $F: \Rmn \rightarrow \Rmn$, initial value $A_0 \in \Rmn$, stepsize $h > 0$, target rank $r > 0$, oversampling parameter $p \geq 0$, number of power iterations $q\geq0$.
\Statex
\State Sketch a Gaussian random matrix $\Omega \in \R^{n \times (r+p)}$.
\State Compute $\Omega^{\dagger} = (\Omega^T \Omega)^{-1} \Omega^T$. \Comment{Solve a small linear system}
\State Compute $B(h)$ by solving the $m \times (r+p)$ differential equation \Comment{B-step}
$$\dt{B}(t) = F(B(t) \Omega^{\dagger}) \Omega, \quad B(0) = A_0 \Omega.$$
\State Extract the basis: $Q = \mathrm{orth}(B(h))$. \Comment{e.g. by pivoted QR}
\For{$i = 1$ to $q$} \Comment{Power iteration loop (optional)}
\State Compute $C(h)$ by solving the $n \times (r+p)$ differential equation
$$\dt{C}(t) = F(Q C(t)^T)^T Q, \quad C(0) = A_0^T Q.$$
\State Extract the basis: $W = \mathrm{orth}(C(h))$.
\State Compute $B(h)$ by solving the $m \times (r+p)$ differential equation
$$\dt{B}(t) = F(B(t) W^T) W, \quad B(0) = A_0 W.$$
\State Extract the basis: $Q = \mathrm{orth}(B(h))$.
\EndFor
\State \textbf{Return} $Q \in \R^{m \times (r+p)}$
\end{algorithmic}
\end{algorithm}

\subsection{Adaptive dynamical rangefinder} \label{sec: adaptive dynamical rangefinder}
In most applications, the rank is not known a priori. 
In a rank adaptive method, we prescribe an error tolerance and let the method choose a suitable rank such that the error $\norm{A(h) - Q_h Q_h^T A(h)}_2$ will be smaller than the desired tolerance.
The so-called adaptive randomized rangefinder \cite[Algorithm 4.2]{halko2011finding} relies on a lemma from~\cite{woolfe2008fast}:
\begin{lemma}[Approximate norm estimate via Gaussian sampling]
Let $M \in \Rmn$, choose an integer $K \geq 1$ and a real number $\alpha > 1$.
Draw an independent family $\{ \omega^{(i)} : i = 1, 2, \ldots, K\}$ of standard Gaussian vectors.
Then,
$$\norm{M}_2 \leq \alpha \sqrt{\frac{2}{\pi}} \cdot \max_{i= 1, \ldots, K} \norm{M \omega^{(i)}}_2$$
with probability at least $1-\alpha^{-K}$.
\end{lemma}
Taking $M = A(h) - QQ^T A(h)$ and $\alpha=10$ leads to the following estimate:
$$\norm{A(h) - QQ^TA(h)}_2 \leq 10 \sqrt{\frac{2}{\pi}} \cdot \max_{i= 1, \ldots, K} \norm{A(h) \omega^{(i)} - QQ^T A(h) \omega^{(i)}}_2$$
with probability at least $1-10^{-K}$.
Recall that the dynamical rangefinder makes the approximation $B(h) \approx A(h) \Omega$ and therefore
\begin{equation*}
 A(h) \omega^{(i)} = A(h) \Omega[:, i] \approx B(h)[:, i] \qquad \text{ where } [:,i] \text{ denotes the $i$th column.}
\end{equation*}
In other words, the columns of $B(h)$ can be used as an accurate error estimate at little extra cost. We define the \textit{adaptive dynamical rangefinder} algorithm as follows:
\begin{enumerate}
  \item Set the block sampling size $K$ according to the prescribed tolerance and probability of failure.
  \item Initial estimation of the range: $Q_h = \mathrm{DynamicalRangefinder}(F, A_0, h, K, p=0, q=0)$.
  \item Draw a new random Gaussian matrix $\Omega \in \R^{n \times K}$.
  \item Solve the $m \times K$ differential equation
    \begin{equation*}
      \left\{
    \begin{aligned}
    &\dt{B}(t) = F(B(t) \Omega^{\dagger}) \Omega, \quad t \in [0,h] \\
    &B(0) = A_0 \Omega. 
    \end{aligned}
      \right.
    \end{equation*}
  \item Compute the error estimate $E = \max_{i = 1, \ldots, K} \norm{B(h)[:, i] - Q_h Q_h^T B(h)[:, i]}_2$.
  \item If the error is too large, augment $Q_h$ with the new basis elements in $B(h)$ and repeat from step 3.
\end{enumerate}
The above procedure is detailed in Algorithm \ref{algo: adaptive dynamical rangefinder}.
\begin{algorithm}
\caption{Adaptive dynamical rangefinder}
\label{algo: adaptive dynamical rangefinder}
\begin{algorithmic}[1]
\Require Vector field $F: \Rmn \rightarrow \Rmn$, initial value $A_0 \in \Rmn$, stepsize $h > 0$, tolerance $\tau > 0$, max. failure probability $\beta > 0$.
\Statex
\State Set $K = - \lceil \log(\beta) / \log(10) \rceil$ and $\varepsilon = \sqrt{\frac{\pi}{2}} \cdot \frac{\tau}{10}$. \Comment{Block size and tolerance}
\State Set $r = K$ and $p=0$.
\State Compute $Q_h = \mathrm{DynamicalRangefinder}(F, A_0, h, r, p, q=0)$.
\State Set $E = \infty$ and $j=1$
\While{$E > \varepsilon$}
\State Sketch a random Gaussian matrix $\Omega \in \R^{n \times K}$ and compute $\Omega^{\dagger} = (\Omega^T \Omega)^{-1} \Omega^T$.
\State Compute $B(h)$ by solving the $m \times K$ differential equation \Comment{B-step}
$$\dt{B}(t) = F(B(t) \Omega^{\dagger }) \Omega, \quad B(0) = A_0 \Omega.$$ 
\State Update the error $$E = \max_{i=1, \ldots, K} \norm{B(h)[:,i] - Q_h Q_h^T B(h)[:,i]}_2.$$
\If{$E > \varepsilon$} \Comment{Augment the basis}
\State Set $j = j+1$ and update the basis: $Q_h = \mathrm{orth}([Q_h, B(h)])$ \Comment{e.g. by a QR update}
\EndIf
\EndWhile
\State \textbf{Return} $Q_h \in \R^{m \times jK}$
\end{algorithmic}
\end{algorithm}

\subsection{Finding the range: a toy experiment}
\begin{figure}
\centering
\includegraphics[width=0.8\textwidth]{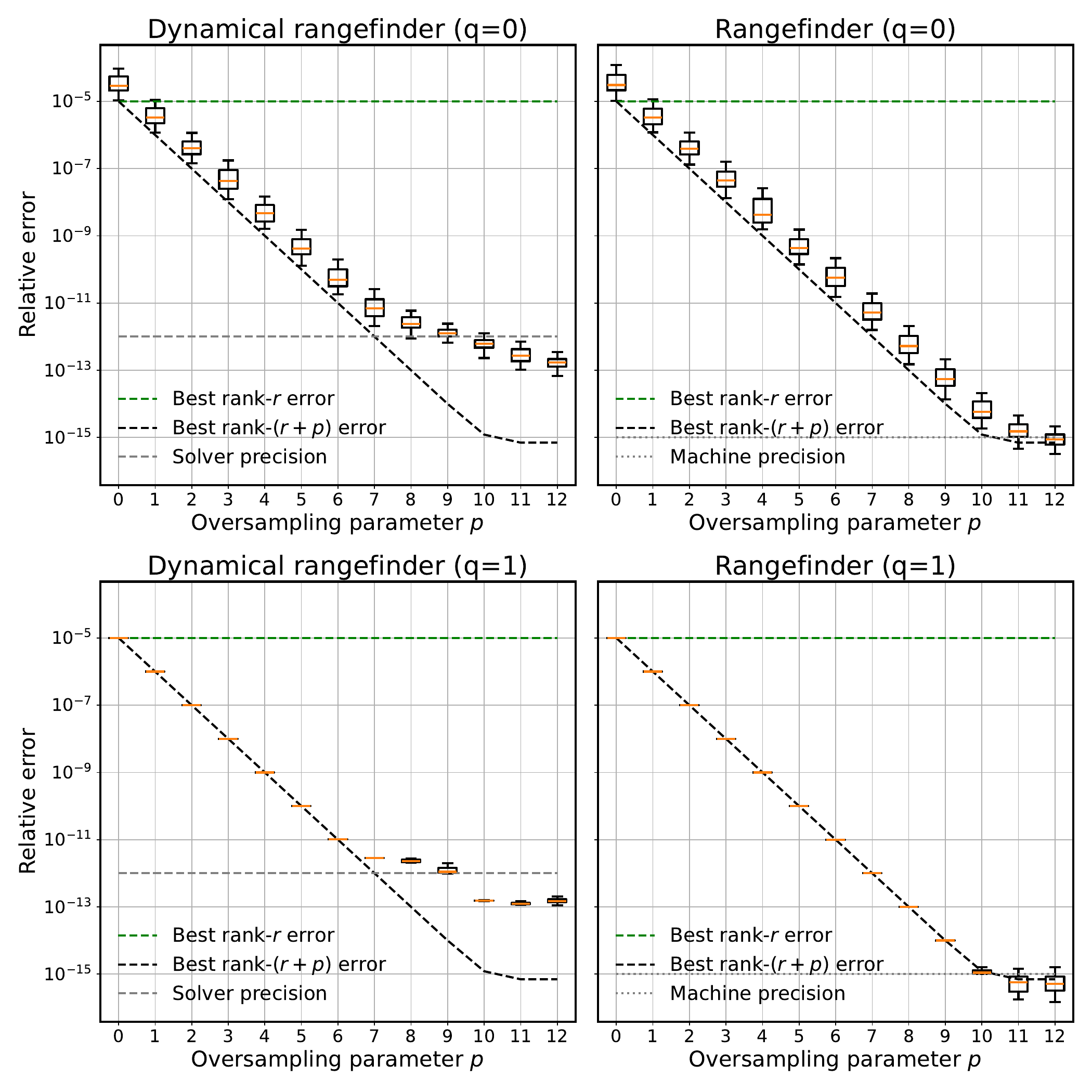}
\caption{Comparison of the dynamical rangefinder algorithm \ref{algo: dynamical rangefinder} with the rangefinder \cite{halko2011finding} applied to the toy problem \ref{eq: toy problem}. The target rank of approximation is $r=5$ and the figure shows the relative error (at time $h=0.1$) made by the methods as we increase the oversampling parameter. The boxes show the median, first and third quartiles and extreme values of the $100$ simulations.}
\label{fig: dynamical rangefinder}
\end{figure}

Let us consider the following differential equation:
\begin{equation} \label{eq: toy problem}
\dt{X}(t) = W_1 X(t) + X(t) + X(t) W_2^T, \quad X(0) = D,
\end{equation}
where $W_1, W_2 \in \R^{n \times n}$ are two random antisymmetric matrices, and $D \in \Rnn$ is a diagonal matrix such that $D_{ii} = 2^{-i}$.
The closed form solution is given by
\begin{equation} \label{eq: toy problem solution}
X(t) = e^{t W_1} \  e^t D \  e^{t W_2^T},
\end{equation}
which we use to compare the idealized and practical methods above defined.  
In the following, we set the parameters to $n=100$ and $h=0.1$. 
On the one hand, we perform the idealized method by applying the rangefinder algorithm to the closed form solution $X(h)$ and we extract the matrix $Q$ with orthonormal columns with a QR decomposition.
On the other hand, we perform the new dynamical rangefinder algorithm \ref{algo: dynamical rangefinder} applied to equation \eqref{eq: toy problem} and call $Q_h$ the resulting matrix with orthonormal columns. We emphasize that with the dynamical rangefinder algorithm, we never solve the full order model but only the reduced differential equation as defined in \eqref{eq: B-step}.
The target rank is $r=5$ and the oversampling parameter $p$ varies from $0$ to $12$.
Due to the random nature of the algorithm, we perform $100$ simulations and show a box-plot showing the relative errors made by each method compared to the reference solution
$$e_1 = \frac{\norm{X(h) - Q Q^T X(h)}_F}{\norm{X(h)}_F}, \quad e_2 = \frac{\norm{X(h) - Q_h Q_h^T X(h)}_F}{\norm{X(h)}_F}.$$
The results are shown in Figure \ref{fig: dynamical rangefinder}. 
As we can see, the two methods have a similar average accuracy, and are close to the rank $(r+p)$ approximation.
In both cases, the variance and median error are improved when a power iteration is performed.
When the oversampling parameter $p$ is large enough, the dynamical rangefinder accuracy is limited by the solver precision (here scipy RK45 with tolerance set to $10^{-12}$). 
This behavior is similar to the limit set by the machine precision for the rangefinder.

\newpage

\section{Time-stepping methods} \label{sec: time stepping mehtods}

Recall again that we are solving the full order model \eqref{eq: full order model} given by the differential equation
$$\dt{A}(t) = F(A(t)), \quad A(0) = A_0 \in \Rmn.$$
In the following, we will assume that the problem can be well-approximated by a low-rank approximation, and in particular we will start from an approximate of the initial value $A_0 \approx Y_0 = U_0 S_0 V_0^T \in \Mr$ where $U_0 \in \Rmr$ and $V_0 \in \Rnr$ have orthonormal columns, and $S_0 \in \Rrr$ is invertible. We could take for example $Y_0 = \Tr(A_0)$, or an approximate of it.
We fix a stepsize $h>0$ and a target rank of approximation $r$, the goal is to find factors such that $A(h) \approx Y_1 = U_1 S_1 V_1^T \in \Mr$ where $U_1$ and $V_1$ have orthonormal columns.
Building on top of the new dynamical rangefinder algorithm, we propose two new time-stepping methods which evolve such a low-rank approximation and which do not involve the projected dynamics~\eqref{eq: dlra model}. We only present one step of each method, which the time-stepping method will repeat to get an approximate solution on the whole time interval. Ultimately, we discuss a simple modification of the two methods that allows to update and automatically choose the rank at each time step.

\subsection{Dynamical randomized SVD} \label{sec: DRSVD}

We start by proposing a new method inspired by the randomized SVD which consists of two steps: (1) a range estimation via the dynamical rangefinder and (2) a post-processing step done on the dynamics projected onto the estimated range. We start by describing an idealized version of the algorithm where we use the exact solution of the full order problem, then we define the practical variant of the algorithm, and numerical experiments are presented in Section~\ref{sec: applications}.

\subsubsection*{An idealized dynamical randomized SVD}

Let us start by motivating the method with an idealized procedure. Consider a target rank $r$ and an oversampling parameter $p$ and let $\mathcal Q_h \in \R^{m \times (r+p)}$ be the matrix with orthonormal columns obtained from the idealized dynamical rangefinder method described in Section~\ref{sec: idealized dynamical rangefinder}. We propose the following idealized procedure:
\begin{itemize}
  \item Solve the differential equation
  \begin{equation} \label{eq: idealized C-step}
    \left\{
    \begin{aligned}
      &\dt{\mathcal C}(t) = F(A(t))^T \mathcal Q_h, \quad t \in [0, h] \\
      &\mathcal C(0) = A_0^T \mathcal Q_h.
    \end{aligned}
    \right.
  \end{equation}
  \item Compute a reduced SVD: $\tilde U_h \Sigma_h V_h^T = \mathrm{svd}(\mathcal C(h)^T)$.
  \item Update the approximation: $A(h) \approx \mathcal A_1 = \mathcal Q_h \tilde  U_h  \Sigma_h  V_h^T$.
\end{itemize}

Additionally, one can truncate the computed solution to the target rank $r$, but is is not strictly necessary here.
A direct consequence of the exact range property stated in Lemma \ref{lemma: exact range property} is the exactness of the above procedure, which we state in the lemma below.
\begin{lemma}[Exactness property] \label{lemma: exactness property}
  Choose a target rank $r$, an oversampling parameter $p$, and let $\mathcal A_1$ be the approximation obtained from the above idealized procedure. 
  If $\mathrm{rank}(A(h)) \leq r+p$, then $\mathcal A_1 = A(h)$.
\end{lemma}
\begin{proof}
  By definition, the solution to the differential equation \eqref{eq: idealized C-step} verifies
  $$\mathcal C(h) = A_0^T \mathcal Q_h + \left[ A(h)^T - A_0^T \right] \mathcal Q_h = A(h)^T Q_h.$$
  Then, we get by construction that $A_1 = \mathcal Q_h \tilde U_h \Sigma_h V_h^T = \mathcal Q_h \mathcal Q_h^T A(h)$ and the conclusion follows with Lemma \ref{lemma: exact range property}, since then $\Range(\mathcal Q_h) = \Range(A(h))$.
\end{proof}
Solving differential equation \eqref{eq: idealized C-step} requires access to the full order model solution, and is therefore not feasible in practice. Instead, we propose to approximate the solution by projecting the differential equation on an approximate range.

\subsubsection*{A practical dynamical randomized SVD}

Because of the time discretization, approximating equation \eqref{eq: idealized C-step} on the estimated range at time $h$ might lead to a poor numerical approximation. We therefore need to augment the approximation space with the information from the initial value. Doing so increases the computational cost of the method but also increase its accuracy. Starting from the initial approximation $A_0 \approx Y_0 = U_0 S_0 V_0^T \in \Mr$, we define one step of the \textit{dynamical randomized SVD} (DRSVD) as follows:
\begin{enumerate}
  \item Estimate the range at time $h$ and augment the space:
  \begin{equation} \label{eq: drsvd range estimation}
    Q_h = \mathrm{DynamicalRangefinder}(F, Y_0, h, r, p, q) \quad \leadsto Q = \mathrm{orth}([U_0, Q_h]).
  \end{equation}
  \item (C-step) Solve the $n \times (2r + p)$ differential equation
  \begin{equation} \label{eq: C-step}
    \left\{
    \begin{aligned}
      &\dt{C}(t) = F(Q C(t)^T)^T Q, \quad t \in [0, h] \\
      &C(0) = Y_0^T Q.
    \end{aligned}
    \right.
  \end{equation}
  \item Compute a reduced SVD of the $(2r+p) \times n$ matrix: $\tilde U_h \Sigma_h V_h^T = \mathrm{svd}(C(h)^T)$, and project the range: $U_h = Q \tilde U_h$.
  \item Truncate the approximated solution to the target rank: $A(h) \approx Y_1^{\mathrm{DRSVD}} = \Tr(U_h \Sigma_h V_h^T)$.
\end{enumerate} 
The truncation is cost-free since the input matrix is already an SVD factorization. 
Moreover, all matrices can be stored by their low-rank factorization, which reduces the algorithm memory footprint.
The most expensive steps are to solve the two differential equations involved, but such steps are typical in dynamical low-rank integrators (see e.g. \cite{lubich2014projector,ceruti2022unconventional,einkemmer2018low}).
The above procedure also verifies the exactness property with an extra assumption on the estimated range.
\begin{lemma}[Exactness property of DRSVD] \label{lemma: exactness property DRSVD}
  Assume that the solution to the full order model \eqref{eq: full order model} at time $h$ verifies $\mathrm{rank}(A(h)) \leq r + p$, and assume that $A_0 = Y_0 = U_0 S_0 V_0^T$. Let $Q$ be as defined in~\eqref{eq: drsvd range estimation}. If for all $t \in [0,h]$,
  $$\Range(A(t)) \subset \Range(Q),$$
  then the dynamical randomized SVD algorithm is exact: $A(h) = Y_1^{\mathrm{DRSVD}}$.
\end{lemma}
\begin{proof}
  Start by defining $L(t) = A(t)^T Q$ and note that $A(t)^T~=~A(t)^T Q Q^T~=~L(t) Q^T$. In particular, it verifies the differential equation
  \begin{equation*}
    \dt{L}(t) = F(A(t))^T Q = F(Q L(t)^T)^T Q, \quad L(0) = A(0)^T Q.
  \end{equation*}
  We therefore have $C(0) = L(0)$ and $\dt{C}(t) = \dt{L}(t)$ for all $t \in [0, h]$, which implies that $C(t) = A(t)^T Q$ and the conclusion follows since $Y_1^{\mathrm{DRSVD}} = \Tr(Q C(h)^T) = \Tr(Q Q^T A(h)) = \Tr(A(h)) = A(h)$.
\end{proof}
The assumption made in the above lemma cannot be checked without having access to the exact solution of the full order model. However, we can already say that the assumption in satisfied in the limit case $(r+p) = \min\{m ,n\}$. In practice, there are two possibilities to increase the probability of satisfying the assumption and get (close to) the exactness property: (i) increase the oversampling parameter $p$ and (ii) reduce the stepsize $h$.
The detailed procedure is summarized in Algorithm \ref{algo: dynamical randomized SVD}.
\begin{algorithm}[H]
\caption{Dynamical randomized SVD (DRSVD)}
\label{algo: dynamical randomized SVD}
\begin{algorithmic}[1]
\Require Vector field $F: \Rmn \rightarrow \Rmn$, initial value $Y_0 = U_0 S_0 V_0^T \in \Mr$, stepsize $h > 0$, target rank $r > 0$, oversampling $p \geq 0$, number of power iterations $q\geq0$.
\Statex
\Statex \textbf{Stage A: Dynamical range estimation}
\State Compute $Q_h = \mathrm{DynamicalRangefinder}(F, Y_0, h, r, p, q)$ \Comment{See Algorithm \ref{algo: dynamical rangefinder}}
\State Augment the basis: $Q = \mathrm{orth}([U_0, Q_h])$ \Comment{With pivoted QR, or else}
\Statex
\Statex \textbf{Stage B: Dynamical post-processing}
\State Compute $C(h)$ by solving the $n \times (2r+p)$ differential equation \Comment{C-step}
$$\dt{C}(t) = F(Q C(t)^T)^T Q, \quad C(0) = Y_0^T Q.$$
\State Compute the reduced SVD of the $n \times (2r+p)$ matrix $C(h)^T = \tilde U \Sigma_1 V_1^T$ and set $U_1 = Q \tilde U$.
\State \textbf{Return} The SVD truncated to the target rank $r$: $A(h) \approx Y_1^{\mathrm{DRSVD}} = \Tr(U_1 \Sigma_1 V_1^T)$.
\end{algorithmic}
\end{algorithm}

\subsection{Dynamical generalized Nyström}

Similarly to the dynamical randomized SVD presented above, we could give an idealized generalized Nyström algorithm based on the exact solution, which would also verify the exactness property stated in Lemma \ref{lemma: exactness property}. 
The ideas being very similar, we skip the discussion of the idealized method and present only the practical algorithm.
Consider a target rank $r>0$, two oversampling parameters $p$, $\ell \geq 0$, and a low-rank approximate of the initial value $A_0 \approx Y_0 = U_0 S_0 V_0^T \in \Mr$. 
We define one step of the \textit{dynamical generalized Nyström} (DGN) algorithm as follows:
\begin{enumerate}
  \item In parallel, estimate the range and corange, and augment the bases:
  \begin{equation} \label{eq: DGN rangefinder step}
    \begin{aligned}
    &Q_h = \mathrm{DynamicalRangefinder}(F, Y_0, h, r, p, q)  &&\leadsto Q = \mathrm{orth}([U_0, Q_h]), \\
    &W_h = \mathrm{DynamicalRangefinder}(F^T, Y_0^T, h, r, p+\ell, q) &&\leadsto W = \mathrm{orth}([V_0, W_h]).
    \end{aligned}
  \end{equation}
  \item (BCD-step) In parallel, solve the three reduced differential equations for $t \in [0,h]$:
  \begin{equation} \label{eq: BCD step}
    \left\{
    \begin{aligned}
    &\dt{B}(t) = F(B(t) W^T) W, \\ 
    &B(0) = A_0 W,
    \end{aligned} 
    \right.
    \qquad
    \left\{
    \begin{aligned}
    &\dt{C}(t) = F(Q C(t)^T)^T Q, \\ 
    &C(0) = A_0^T Q,
    \end{aligned} 
    \right.
    \qquad
    \left\{
    \begin{aligned}
    &\dt{D}(t) = Q^T F(Q D(t) W^T) W, \\ 
    &D(0) = Q^T A_0 W.
    \end{aligned} 
    \right.
    \end{equation}
  \item Truncate and assemble the solution:
  \begin{equation} \label{eq: DGN truncate and assemble}
    A(h) \approx Y_1^{\mathrm{DGN}} = B(h) \ \Tr(D(h))^{\dagger} \ C(h)^T.
  \end{equation}
\end{enumerate}
As for the generalized Nyström algorithm given in \eqref{eq: orthogonal generalized Nyström}, the order in which the solution is assembled is important for the numerical stability of the method. We first truncate and only then efficiently compute the pseudoinverse of the truncated matrix.
Similarly to the dynamical randomized SVD, the dynamical generalized Nyström algorithm verifies the following exactness property when we assume that the augmented range and corange contain respectively the range and corange of the exact solution.
\begin{lemma}[Exactness property of DGN]
  Assume that the solution to the full order model \eqref{eq: full order model} at time $h$ verifies $\mathrm{rank}(A(h)) \leq r + p$, and assume that $A_0 = Y_0 = U_0 S_0 V_0^T$.
  Let $Q$ and $W$ be the matrices as defined in \eqref{eq: DGN rangefinder step}.
  If for all $t \in [0, h]$,
  $$\Range(A(t)) \subset \Range(Q) \quad \text{and} \quad \Range(A(t)^T) \subset \Range(W),$$
  then the dynamical generalized Nyström algorithm is exact: $A(h) = Y_1^{\mathrm{DGN}}$.
\end{lemma}
\begin{proof}
  The argument made in the proof of Lemma \ref{lemma: exactness property DRSVD} holds verbatim on the three differential equations~\eqref{eq: BCD step}, so we have the three quantities:
$$ B(h) = A(h) W, \quad C(h)^T = Q A(h)^T, \quad D(h) = Q^T A(h) W.$$
The conclusion follows from the definitions and noticing that
$$Y_1^{\mathrm{DGN}} = B(h) \ \Tr(D(h))^{\dagger} \ C(h)^T = A(h) W (Q^T A(h) W)^{\dagger} Q A(h)^T = A(h),$$
thanks again to the range inclusions and the assumption on the rank of $A(h)$.
\end{proof}
The detailed algorithm is summarized in Algorithm \ref{algo: dynamical generalized Nystrom}, and we refer to Section~\ref{sec: applications} for applications.
\begin{algorithm}[H]
\caption{Dynamical generalized Nyström (DGN)}
\label{algo: dynamical generalized Nystrom}
\begin{algorithmic}[1]
\Require Vector field $F: \Rmn \rightarrow \Rmn$, initial value $Y_0 = U_0 S_0 V_0^T \in \Mr$, stepsize $h > 0$, target rank $r>0$, oversampling parameters $p \geq 0$ and $\ell \geq 0$, number of power iterations $q \geq 0$.
\Statex
\Statex \textbf{Stage A: Dynamical range and corange estimation}
\State Compute (in parallel): \Comment{See Algorithm \ref{algo: dynamical rangefinder}.}
\begin{align*}
&Q_h = \mathrm{DynamicalRangefinder}(F, Y_0, h, r, p, q) \\
&W_h = \mathrm{DynamicalRangefinder}(F^T, Y_0^T, h, r, p+\ell, q) 
\end{align*} 
\State (Optional) Augment the two bases: $Q = \mathrm{orth}([U_0, Q_h]), W = \mathrm{orth}([V_0, W_h])$.
\Statex
\Statex \textbf{Stage B: Solve three small differential equations in parallel} \Comment{BCD-step}
\State Compute $B(h)$ by solving the $m \times (2r+p+\ell)$ differential equation
$$\dt{B}(t) = F(B(t) W^T) W, \quad B(0) = Y_0 W, \qquad t \in [0,h].$$
\State Compute $C(h)$ by solving the $n \times (2r+p)$ differential equation
$$\dt{C}(t) = F(Q C(t)^T)^T Q, \quad C(0) = Y_0^T Q, \qquad t \in [0,h].$$
\State Compute $D(h)$ by solving the $(2r+p) \times (2r+p+\ell)$ differential equation
$$\dt{D}(t) = Q^T F(Q D(t) W^T) W, \quad D(0) = Q^T Y_0 W, \qquad t \in [0,h].$$
\Statex \textbf{Stage C: Truncate and assemble according to Nyström's formula}
\State Compute a reduced SVD of the $(2r+p) \times (2r+p+\ell)$ matrix: $\tilde U \Sigma \tilde V^T = \mathrm{svd}(D(h))$.
\State Truncate it to the target rank: $\mathcal T_r(D(h)) = \tilde U_r \Sigma_r \tilde V^T_r$.
\State Compute two QR decompositions and assemble:
$$U_1R_1 = B(h) \tilde V_r, \quad V_1 R_2 = C(h) \tilde U_r \qquad \leadsto S_1 = R_1 \Sigma_r^{-1} R_2^T.$$
\State \textbf{Return} $A(h) \approx Y_1^{\mathrm{DGN}} = U_1 S_1 V_1^T$.
\end{algorithmic}
\end{algorithm}

\newpage

\subsection{Rank-adaptive time-stepping methods}

In the two new methods above presented, a fixed target rank is prescribed, which assumes that we have a guess on the numerical rank of the solution. In practice, we do not have always such a guess and we typically need to change the rank at each time step. Unlike the projector-splitting integrator, the augmented space together with the oversampling parameter $p$ already allow some freedom in the truncation step so that the rank could evolve from one step to the next. However, the rank growth would be limited to $2r+p$, which might not be enough to catch fast changing solutions.
Instead, we propose a more systematic approach based on the adaptive dynamical rangefinder algorithm proposed in Section~\ref{sec: adaptive dynamical rangefinder}, which allows to choose automatically the rank.

Consider a prescribed error tolerance $\tau>0$ and a certain probability of failure $\beta>0$. Let us denote the tolerance-based truncated SVD by $\mathcal T_{\tau}$ and defined such that 
$$\norm{X - \mathcal T_{\tau}(X)}_2 \leq \tau \quad \text{ for any matrix } X \in \Rmn \text{ and tolerance } \tau \in \R,$$
and which can be computed by truncating the singular values of the matrix $X$. 
Starting from the initial low-rank approximate $Y_0 = \mathcal T_{\tau}(A_0) = U_0 S_0 V_0^T$, we define one step of the \textit{adaptive dynamical randomized SVD} (ADRSVD) algorithm as follows:
\begin{enumerate}
  \item Estimate the range according to $\tau$ and $\beta$, and augment the basis:
  $$Q_h = \mathrm{AdaptiveDynamicalRangefinder}(F, Y_0, h, \tau, \beta) \quad \leadsto Q = \mathrm{orth}([U_0, Q_h]).$$
  \item (C-step) Solve the reduced differential equation \eqref{eq: C-step} to get $C(h)$.
  \item Compute a reduced SVD: $\tilde U_h \Sigma_h V_h^T = \mathrm{svd}(C(h)^T)$, and project the range: $U_h = Q \tilde U_h$.
  \item Truncate the approximated solution to the prescribed tolerance: $A(h) \approx Y_1^{\mathrm{ADRSVD}} = \mathcal T_{\tau} (U_h \Sigma_h V_h^T)$.
\end{enumerate}
The definition of the \textit{adaptive dynamical generalized Nyström} (ADGN) is similar: 
\begin{enumerate}
  \item In parallel, estimate the range and corange according to $\tau$ and $\beta$, and augment the bases:
    \begin{equation} \label{eq: SDGN rangefinder step}
    \begin{aligned}
    &Q_h = \mathrm{AdaptiveDynamicalRangefinder}(F, Y_0, h, \tau, \beta)  &&\leadsto Q = \mathrm{orth}([U_0, Q_h]), \\
    &W_h = \mathrm{AdaptiveDynamicalRangefinder}(F^T, Y_0^T, h, \tau, \beta) &&\leadsto W = \mathrm{orth}([V_0, W_h]).
    \end{aligned}
  \end{equation}
  \item (BCD-step) Solve the three reduced differential equations \eqref{eq: BCD step} to get $B(h)$, $C(h)$, and $D(h)$.
  \item Truncate up to the prescribed tolerance and assemble the solution:
  \begin{equation} \label{eq: ADGN truncate and assemble}
    A(h) \approx Y_1^{\mathrm{ADGN}} = B(h) \ \mathcal T_{\tau}(D(h))^{\dagger} \ C(h)^T.
  \end{equation}
\end{enumerate}
The implementation of the two adaptive methods above defined is very similar to the implementation given in Algorithm~\ref{algo: dynamical randomized SVD} and Algorithm~\ref{algo: dynamical generalized Nystrom}, and we refer to Section~\ref{sec: application Allen Cahn adaptive} for an application of these two adaptive methods on the Allen--Cahn equation.

\newpage

\section{Applications} \label{sec: applications}

We now apply the new time-stepping methods to partial differential equations and we compare the performance (cost and accuracy) with state-of-the-art DLRA techniques.
The implementation is done in Python, and the code is fully available on \href{https://github.com/BenjaminCarrel/randomized-dynamical-low-rank}{GitHub}\footnote{The URL is: https://github.com/BenjaminCarrel/randomized-dynamical-low-rank} for reproducibility.
All experiments have been executed on a MacBook Pro with an M1 processor and 16GB of RAM.

\subsection{A stiff differential equation with Lyapunov structure}

Stiff differential equations are numerically challenging to solve as they impose a strong restriction on the stepsize.
Methods like the projector-splitting integrator \cite{lubich2014projector} and explicit (low-rank) Runge--Kutta methods \cite{kieri2019projection,lam2024randomized} are not suited for such problems. 
We consider here the heat equation with a constant source, given by the partial differential equation
\begin{equation} \label{eq: heat equation}
u(x, y; t) = \Delta u(x, y; t) + s(x, y), \quad u(x, y; 0) = u_0(x, y),
\end{equation}
with Dirichlet boundary conditions and where $(x,y) \in \mathcal D$ and $t \in [0, T]$. 
Here, the source $s(x,y)$ is constant in time.
A spatial discretization via finite differences leads to the Lyapunov differential equation \cite{uschmajew2020geometric} defined by
\begin{equation} \label{eq: lyapunov differential equation}
\dt{A}(t) = L A(t) + A(t) L + \alpha \frac{C}{\norm{C}_F}, \quad A(0) = A_0,
\end{equation}
where $A(t) \in \Rnn$, $\alpha \geq 0$ and $t \in [0, T]$. 
For our numerical experiments, we use $\Omega = [-\pi, \pi]^2$ discretized with $n=256$ points. 
The discrete Laplacian operator is $L = (n-1)^2 \cdot \mathrm{diag}(1, -2, 1) \in \R^{n \times n}$. 
The construction of the initial value and the source is similar to \cite{uschmajew2020geometric,lam2024randomized}.
Since Assumption~\ref{ass: closeness} is not satisfied at the initial value, we propagate the solution from time $0$ to $10^{-4}$ and start our numerical experiments from there. We point out that Assumption~\ref{ass: closeness} is only required for DLRA techniques, and we do not require it for the new methods. The new methods do \textbf{not} require this extra propagation of the initial value.

In Table \ref{table: lyapunov experiment}, we compare the accuracy of the methods applied to equation~\eqref{eq: lyapunov differential equation} with a target rank $r=5$.
When the stepsize is large (here $h=0.1$), existing DLRA methods perform poorly because the problem is stiff.
The new methods perform well in comparison with a good accuracy and low variance despite the randomness. The table shows the importance of the power iteration for the DRSVD algorithm. Note that the DGN algorithm with $q=1$ power iteration is able to recover the reference best rank approximation.

\begin{figure}[!ht]
\centering
\includegraphics[width=0.8\textwidth]{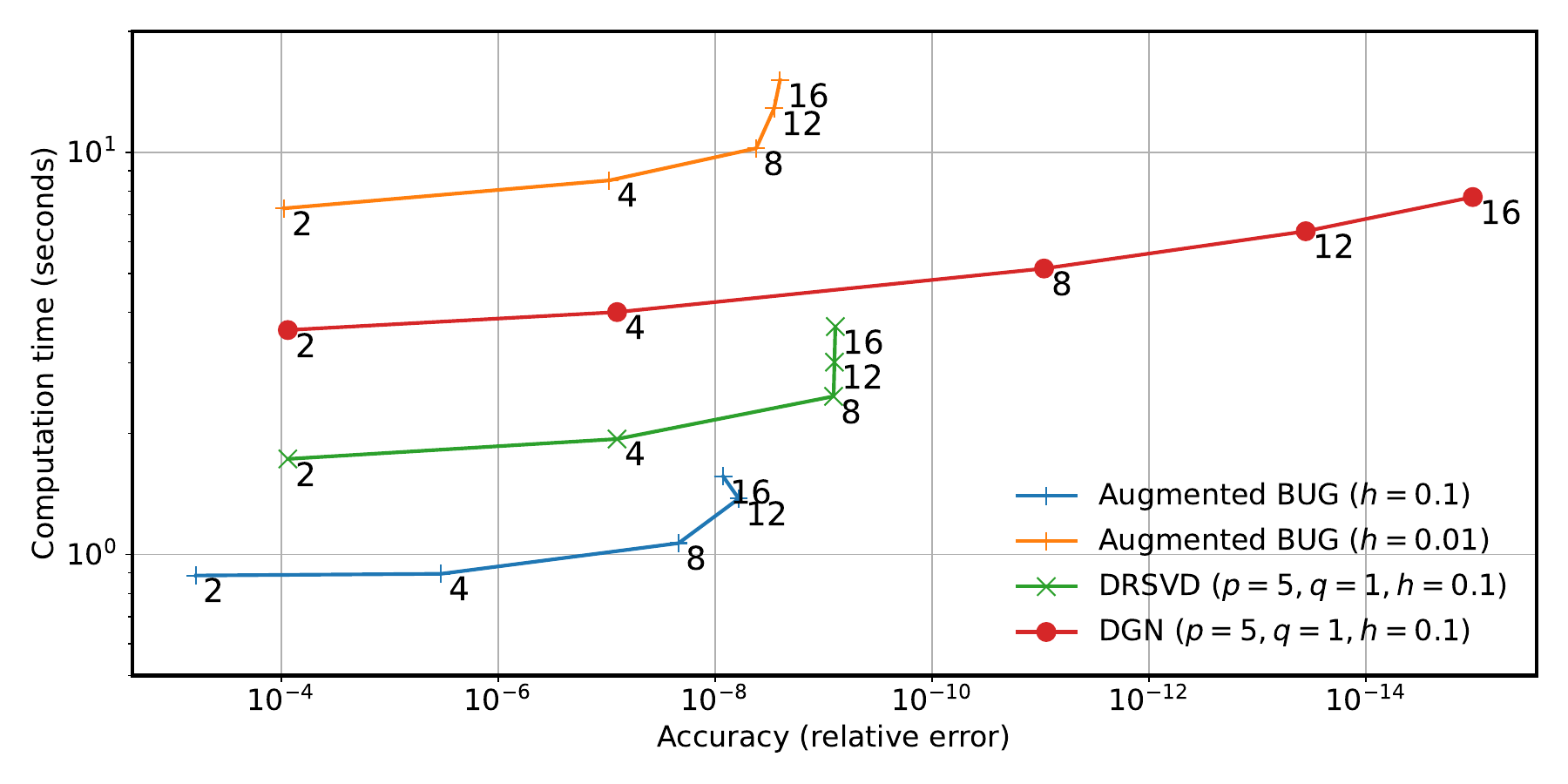}
\caption{Accuracy and computation time (in seconds) of the methods applied to the Lyapunov equation~\eqref{eq: lyapunov differential equation}. The integers next to the markers correspond to the target rank. The reference solver is a second-order exponential Runge--Kutta method with stepsize~$\delta t = 10^{-4}$, and the substeps of the low-rank methods were solved with this reference solver.}
\label{fig: lyapunov performance}
\end{figure}

In Figure~\ref{fig: lyapunov performance}, we report the computational time versus accuracy for the augmented BUG integrator with two different step sizes and for the proposed methods. For $h=0.1$, the augmented BUG is fastest but fails to converge, as the relative error stagnates above our expectations (for example by comparing the errors obtained with a target rank $r=2$). Reducing the stepsize to $h=0.01$ restores convergence, but at the expense of substantially increased computational cost. By contrast, the proposed methods (DRSVD and DGN) achieve significantly higher accuracy at lower cost thanks to a convergence achieved with a larger stepsize. In particular, DGN consistently reaches the highest attainable accuracy (down to $10^{-15}$) when the target rank is sufficiently large, thereby outperforming the augmented BUG integrator in both efficiency and robustness to get convergence.

\begin{table}[H]
\begin{subtable}{\textwidth}
  \centering
  \caption{Existing DLRA methods with a rank of approximation $r=5$, and the reference solution truncated to the same rank.}
  \begin{tabular}{lc}
    \toprule
  Method & Relative error after one step \\
  \cmidrule(lr){1-1}\cmidrule(lr){2-2}
  Projector-splitting \cite{lubich2014projector} & (diverged) \\
  Projected RK1 \cite{kieri2019projection} & 1.49e-01 \\
  Randomized low-rank RK1 \cite{lam2024randomized} & 1.49e-01 \\
  Projected exponential RK1 \cite{carrel2023projected} & 1.01e-04 \\
  BUG integrator \cite{ceruti2022unconventional} & 3.37e-05 \\
  Augmented BUG integrator \cite{ceruti2022rank} & 1.04e-06 \\
  Reference truncated to the target rank ($r=5$) & 4.50e-09 \\
  \end{tabular}
\end{subtable}
\begin{subtable}{\textwidth}
  \centering
  \begin{tabular}{ccccc}
  \toprule
  Method & $p=0$ & $p=2$ & $p=5$ & $p=10$ \\
  \cmidrule(lr){1-1}\cmidrule(lr){2-2}\cmidrule(lr){3-3} \cmidrule(lr){4-4}\cmidrule(lr){5-5}
  DRSVD ($q=0$)  & 3.11e-04 &  1.93e-04 & 1.29e-04 & 8.29e-05 \\
  (Quartiles) &[2.83e-04, 3.31e-04] &[1.84e-04, 1.94e-04] &[1.25e-04, 1.34e-04] &[7.98e-05, 8.69e-05] \\
  DRSVD ($q=1$) & 3.25e-08 & 6.94e-09 & 6.08e-09 & 4.50e-09 \\
  (Quartiles) &[3.25e-08, 3.26e-08] &[6.41e-09, 7.16e-09] &[4.71e-09, 7.22e-09] &[4.50e-09, 4.50e-09] \\
  DGN ($q=0$) & 5.19e-09 & 4.66e-09 & 4.54e-09 & 4.51e-09 \\
  (Quartiles) &[5.01e-09, 5.32e-09] &[4.65e-09, 4.67e-09] &[4.53e-09, 4.55e-09] &[4.51e-09, 4.52e-09] \\
  DGN ($q=1$) & 4.50e-09 & 4.50e-09 & 4.50e-09 & 4.50e-09 \\
  (Quartiles) &[4.50e-09, 4.50e-09] &[4.50e-09, 4.50e-09] &[4.50e-09, 4.50e-09] &[4.50e-09, 4.50e-09] \\
  \bottomrule
  \end{tabular}
\caption{Median, first and third quartiles $([25\%, 75\%])$ relative error ($30$ simulations) made by the new methods: dynamical randomized SVD (DRSVD) and dynamical generalized Nyström (DGN) with oversampling parameter $p$ and $q$ power iterations. For DGN, the second oversampling parameter is $\ell=0$.}
\end{subtable}
\caption{Relative error on the stiff matrix differential equation~\eqref{eq: lyapunov differential equation} after one step with $h=0.1$. The reference solver is a second-order exponential Runge--Kutta method with stepsize~$\delta t = 10^{-4}$, and the substeps were solved with the reference solver.}
\label{table: lyapunov experiment}
\end{table}

\subsection{Reaction-diffusion with the Allen--Cahn equation} \label{sec: application Allen Cahn adaptive}

Another example of a stiff differential equation comes from the process of phase separation in multi-component alloy systems, described by the Allen--Cahn equation \cite{allen1972ground,allen1973correction}. In its simplest form, it is given by
\begin{equation} \label{eq: Allen--Cahn continuous equation}
\frac{\partial f}{\partial t} = \varepsilon \Delta f + f - f^3,
\end{equation}
where $\Delta f$ is the diffusion term and $f - f^3$ is the reaction term.
Dynamical low-rank techniques tailored for stiff differential equations have been successfully applied to this problem, see \cite{rodgers2023implicit,carrel2023projected}. 
The stiffness can be controlled with the small parameter $\varepsilon$, and here we take $\varepsilon = 0.01$ to get a moderately stiff problem.
We consider a setting similar to~\cite{rodgers2023implicit} with initial value
$$f_0(x,y) = \frac{2 \ e^{- \tan^2(x)} \sin(x) \sin(y) }{1 + e^{|\csc(-x/2)|} + e^{|\csc(-y/2)|}},$$
the domain is $\mathcal D = [0, 2\pi]^2$ and the spatial mesh contains $128 \times 128$ points with periodic boundary conditions.
After finite differences, we obtain the matrix differential equation
\begin{equation} \label{eq: Allen--Cahn discrete equation}
\dt{X}(t) = A X(t) + X(t) A + X(t) - X(t)^{*3},
\end{equation}
where the power is taken element-wise (Hadamard product). 
We consider the time interval $[0, T] = [0, 10]$, and most of the reaction occurs between $t=3$ and $t=7$. 
At the end of the time interval, the solution has almost reached the steady-state.
At any given time, the singular values of the solution decay exponentially fast. 
 
\begin{figure}[ht!]
  \centering
\includegraphics[width=0.9\textwidth]{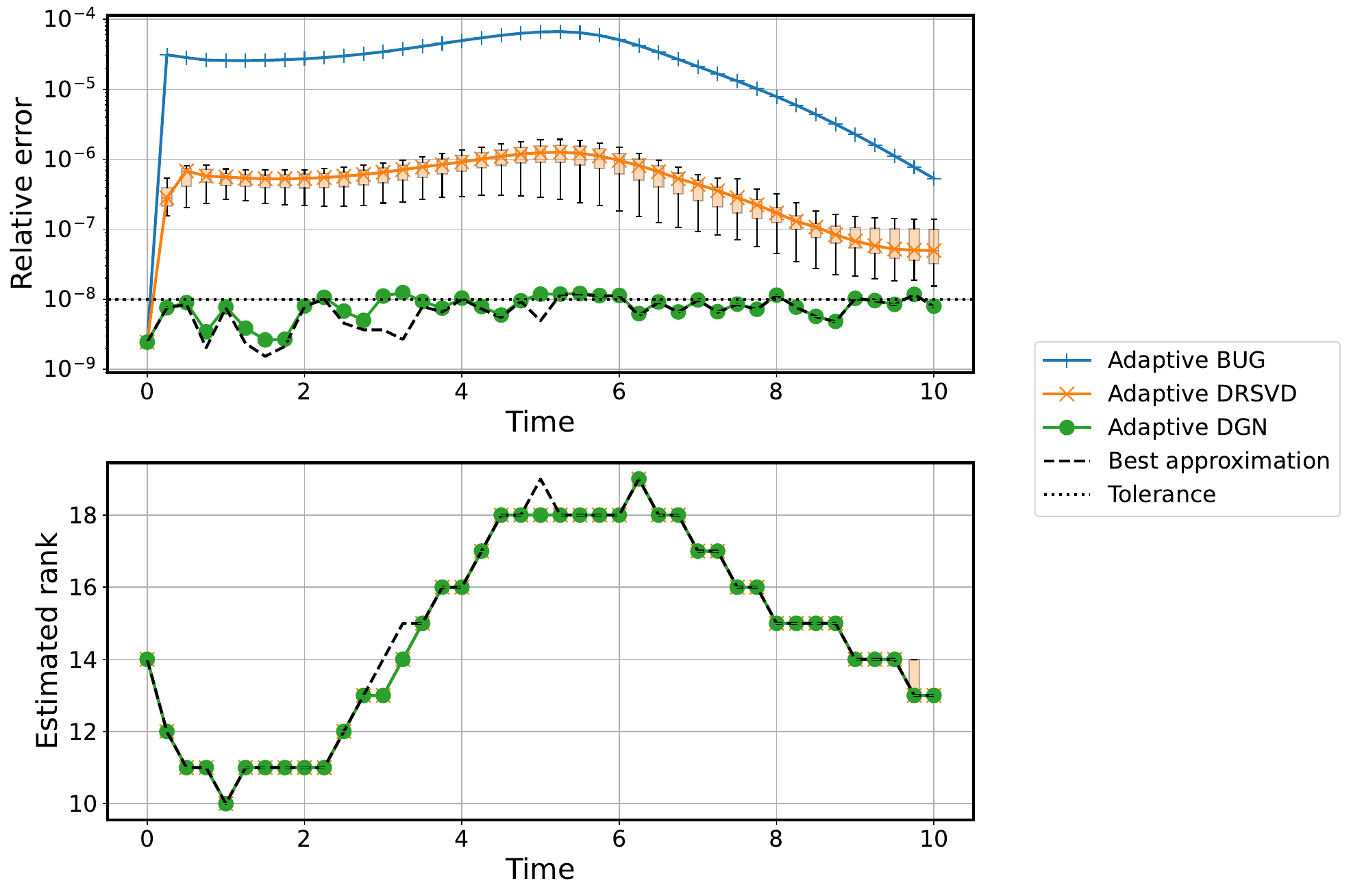}
\caption{Rank-adaptive methods applied to the Allen--Cahn equation \eqref{eq: Allen--Cahn discrete equation} with stepsize $h = 0.25$ and $30$ simulations. The prescribed tolerance is $\tau = 10^{-8}$. The reference is computed with a second order exponential Runge--Kutta method with time step $\delta t = 0.005$, and the substeps in the low-rank methods are also solved with the reference solver.}
\label{fig: allen-cahn experiments}
\end{figure}

The reaction component leads to a non-trivial evolution of the singular values, it is therefore a good problem for testing our new rank-adaptive methods.
In Figure \ref{fig: allen-cahn experiments}, we compare the new time-stepping methods with the rank-adaptive BUG integrator \cite{ceruti2022rank}. We perform $30$ simulations of the random methods and show, at each time step, a boxplot showing the median and quartiles of the relative error. For a stepsize $h=0.5$, the error due to the time discretization is large for the adaptive BUG and adaptive DRSVD. Even though the variance of adaptive DRSVD is large, the method is always more accurate than the adaptive BUG on this example. On the other hand, the adaptive DGN has a tiny variance and a relative error close to the prescribed tolerance. Despite the differences in the relative error, all three methods estimate a rank of approximation close to the reference.

\subsection{Stochastic problems with the Burgers equation}

In uncertainty quantification, the discretization of stochastic partial differential equations leads to large matrix differential equations. 
Let us consider here the one-dimensional Burgers equation
\begin{equation} \label{eq: burgers equation}
\partial_t v + v  \partial_x v  = \nu \partial_{xx} v , \quad x \in \mathcal D, \quad t \in [0, T],
\end{equation}
with Dirichlet boundary conditions.
As done in \cite{naderi2024cur}, we sample the problem on several random initial conditions defined by
\begin{equation} \label{eq: burgers initial condition}
v(x, 0; \xi) = 0.5 \sin(2 \pi x) \left( e^{\cos(2 \pi x)} - 1.5 \right) + \sigma_x \sum_{i=1}^d \sqrt{\lambda_{x_i}} \psi_i(x) \xi_i, \quad \xi_i \sim \mathcal N(0, 1).
\end{equation}
More specifically, we draw $d$ random vectors $\xi_i \in \R^s$ that are normally distributed, and the eigenvalues $\lambda_i$ and eigenvectors $\psi_i(x)$ are extracted from the eigenvalue decomposition of the kernel matrix 
$$K_{ij} = \exp(-(x_i - x_j)^2/2).$$ 
The problem is discretized by associating each column to a random sample. Hence, taking the space $\mathcal D = [0, 1]$ and using central finite differences leads to the matrix differential equation
\begin{equation} \label{eq: matrix burgers equation}
\dt{A}(t) = L A(t) - A(t) * (D A(t)),
\end{equation}
where the unknown matrix is $A(t) \in \R^{n \times s}$ and the operation $*$ is again the element-wise (Hadamard) product. 
The discrete differential operators are 
$$L = (n-1)^2 \cdot \mathrm{diag}(1, -2, 1) \in \R^{n \times n}, \quad D = (n-1)/2 \cdot \mathrm{diag}(-1, 0, 1) \in \R^{n \times n}.$$

\begin{figure}[ht!]
\centering
\includegraphics[width=0.8\textwidth]{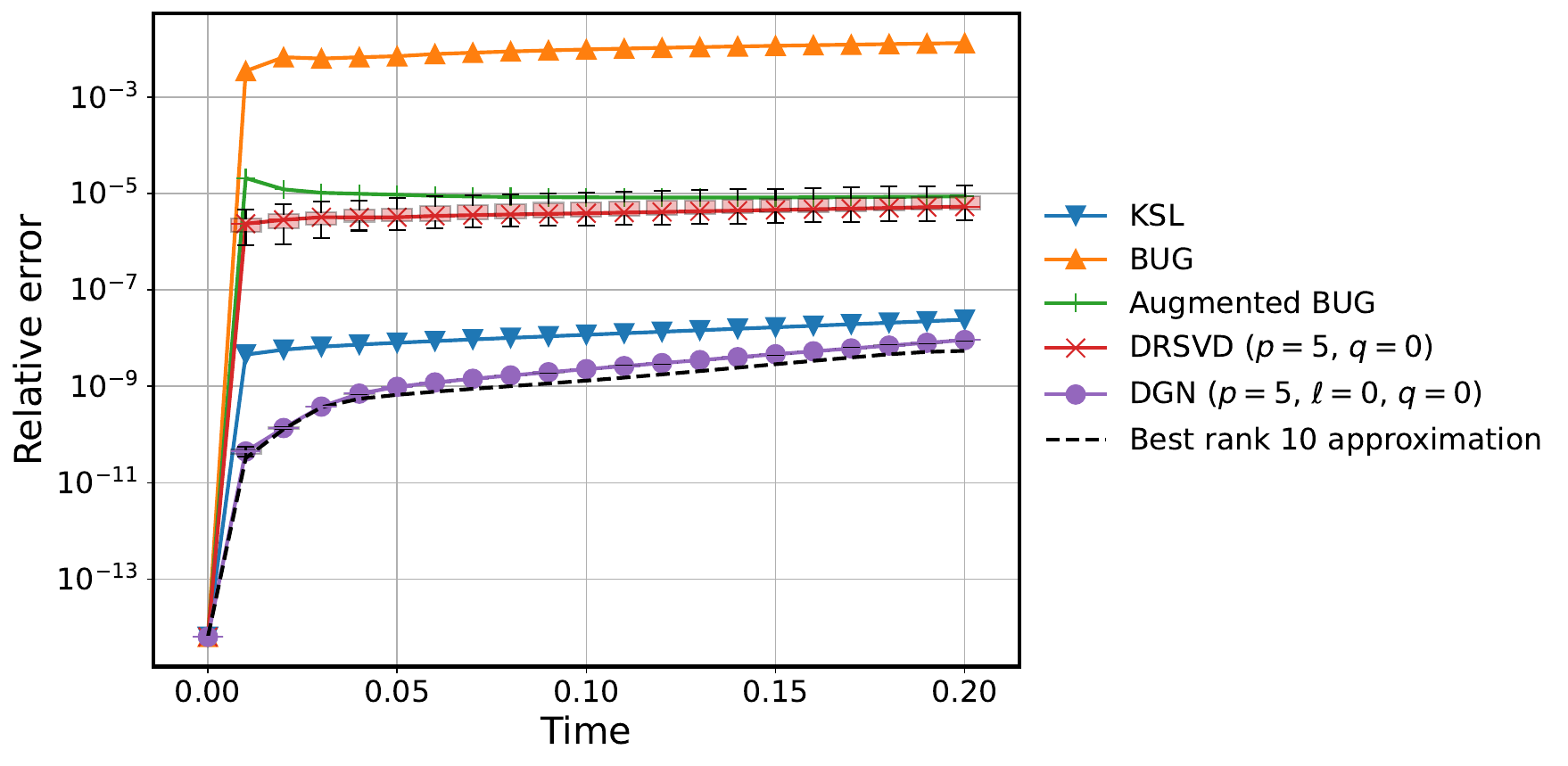}
\caption{Relative error over time of several methods applied to the stochastic Burgers equation \eqref{eq: matrix burgers equation} with stepsize $h=0.01$ and $30$ simulations. The reference solver is scipy RK45 method with absolute and relative tolerance set to $10^{-12}$, and the substeps of the low-rank methods are solved with this reference solver.}
\label{fig: burgers experiment}
\end{figure}

In our experiments, we set the parameters to $\nu = 0.01,$ $\sigma_x = 0.001$ and $d=4$. 
Then, we take $n=256$ spatial points, $s=64$ random samples, and solve the problem until $T=0.2$.
For those parameters, the differential equation is not stiff. 
In Figure \ref{fig: burgers experiment}, we compare dynamical low-rank techniques on the stochastic Burgers equation~\eqref{eq: matrix burgers equation} with a target rank of approximation $r=10$. Since the differential equation is not stiff, the projector-splitting works particularly well and has a small hidden constant compared to the BUG integrators. On this problem, the dynamical randomized SVD has a performance similar to the augmented BUG, and the dynamical generalized Nyström has the best accuracy with a relative error close to the best rank approximation.

\subsection{Kinetic equations with Vlasov-Poisson}

The following Vlasov-Poisson equation~\cite{einkemmer2018low} describes the evolution of electrons in a collisionless plasma, in the electrostatic regime with constant background ion density:
\begin{equation} \label{eq: Vlasov-Poisson}
\begin{aligned}
&\partial_t f(t, x, v) + v \cdot \nabla f(t, x, v) - E(f)(t, x) \cdot \nabla_v f(t, x, v) = 0, \\
&\nabla \cdot E(f)(t, x) = - \int f(t, x, v) dv + 1, \qquad \nabla \times E(f)(t, x) = 0,
\end{aligned}
\end{equation}
where $f(t,x,v)$ is the particle-density function with $x \in \mathcal D_x \subset \R^d$ and $v \in \mathcal D_v \subset \R^d$.
One aspect of equation \eqref{eq: Vlasov-Poisson} is its large memory footprint since it requires to store $O(n^{2d})$ floating point numbers, where $n$ is the number of grid points. 
To overcome the memory limitations, low-rank techniques have been introduced in~\cite{einkemmer2018low} and we refer to references therein for state-of-the-art full rank solvers. The research on low-rank techniques for the Vlasov equation is still in development, see e.g.~\cite{guo2022low,cassini2022efficient,guo2024conservative,uschmajew2024dynamical}. 
Note that in the literature, the discretization is typically done after deriving equations for the low-rank factors of the function $f(t,x,v)$ itself.
In this paper, we first perform a spatial discretization via finite differences and get the following matrix differential equation:
\begin{equation} \label{eq: discretised Vlasov-Poisson}
\dt{A}(t) = - D_x A(t) V - E(t, A(t)) A(t) D_v,
\end{equation}
where $D_x$ and $D_v$ are (high-order) centered finite difference matrices, $V$ and $E(t, A(t))$ are two diagonal matrices with coefficients derived from $\eqref{eq: Vlasov-Poisson}$. 
In the following experiments, the reference solver for solving equation~\eqref{eq: discretised Vlasov-Poisson} is scipy RK45 solver with tolerance~$10^{-8}$.

\subsubsection*{Linear Landau damping}

The Landau damping is a phenomenon that occurs in plasma and which can be modeled with equation~\eqref{eq: Vlasov-Poisson}.
It is a damping of the electric field due to the interaction of the particles with the field. 
We consider the same setup as provided in~\cite{einkemmer2018low}, that is, the domain is $\Omega = (0, 4 \pi) \times (-6, 6)$ with periodic boundary conditions and initial value
\begin{equation} \label{eq: linear landau damping}
f(0, x, v) = \tfrac{1}{\sqrt{2 \pi}} \exp \left( - \tfrac{v^2}{2} \right) \left( 1 + \alpha \cos(kx) \right),
\end{equation}
where $\alpha = 10^{-2}$ is the amplitude of the perturbation and $k = 0.5$ is the wave number.
For the discretization, we used $n_x = 64$ and $n_v = 256$ grid points in the $x$ and $v$ directions, respectively.
For this specific problem, it is known that the theoretical rate of decay for the electric energy is $\gamma \approx -0.153$.

\begin{figure}[ht!]
\centering
\includegraphics[width=0.8\textwidth]{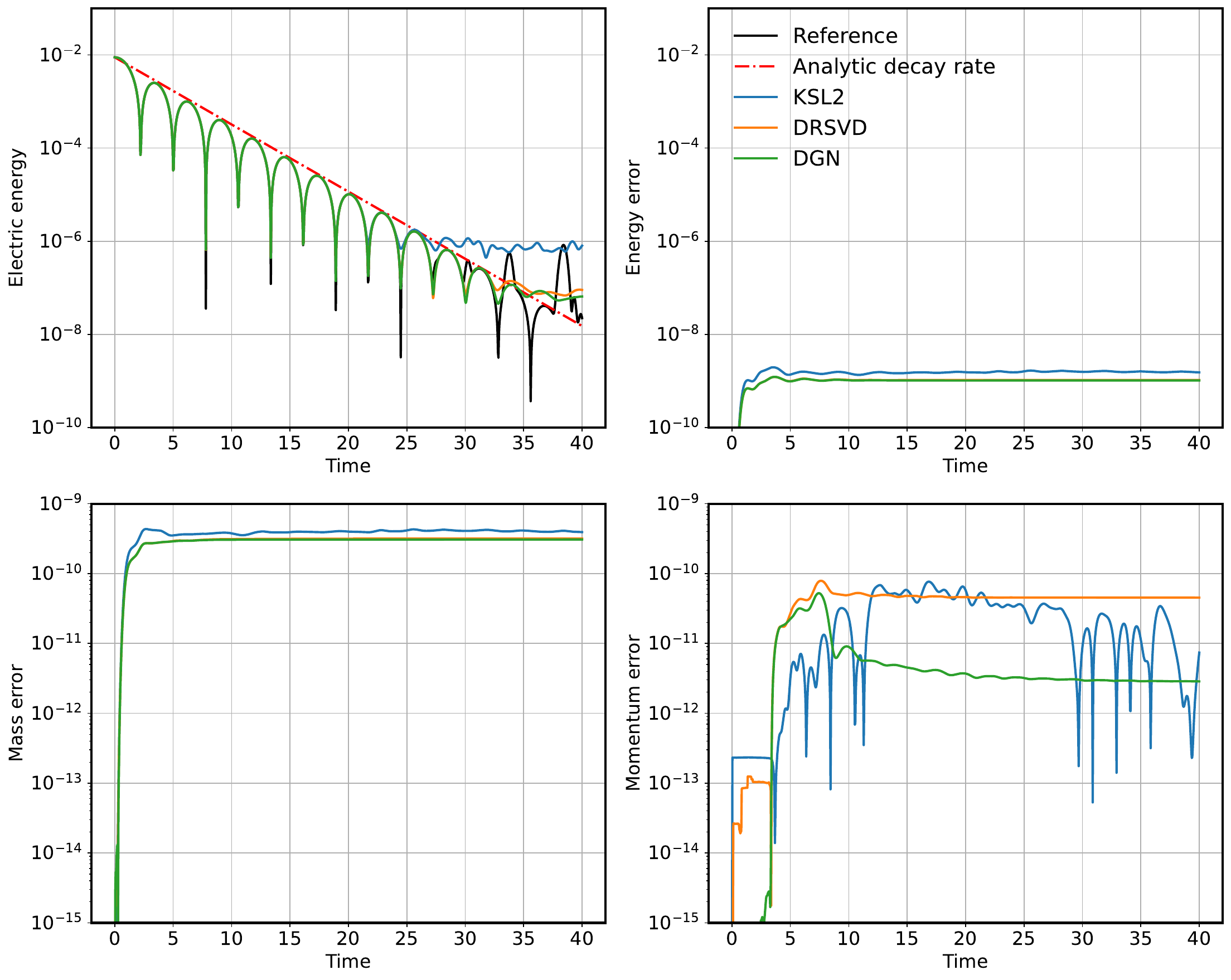}
\caption{Numerical simulations of the linear Landau damping \eqref{eq: linear landau damping} with stepsize $h=0.04$. The projector-splitting is performed with a Strang splitting ($\mathrm{KSL2}$) and the substeps are solved with the reference solver. The DRSVD is performed with oversampling parameter $p=5$ and $q=1$ power iteration. The DGN is performed with oversampling parameters $p=5$, $\ell=0$ and $q=1$ power iteration.}
\label{fig: landau damping}
\end{figure}

In Figure \ref{fig: landau damping}, we apply the new methods to the Landau damping problem until the final time $T=40$. The stepsize is $h=0.04$, and the substeps are solved with the reference solver. The projector-splitting (KSL2) is also given as a reference low-rank technique. The new methods accurately solve the problem and recover the theoretical rate until the tolerance of the reference solver is reached. At the end of the time interval, the modelling error induced by the new techniques is smaller than the modelling error induced by the projector-splitting integrator.

\subsubsection*{Two-stream instability}

The two-stream instability~\cite{bittencourt2013fundamentals} is a well-known instability in plasma physics which is also modeled by equation~\eqref{eq: Vlasov-Poisson}.
It is characterized by two counter-propagating beams of particles, which are perturbed by a small perturbation when $\alpha \neq 0$ in the initial value given below. 
The perturbation in the electric field grows exponentially, leading to the formation of structures in the plasma. 
Ultimately, the problem reaches the non-linear regime characterized by a nearly constant electric energy and filamentation of the phase space. 
While the linear regime admits good low-rank approximations, the non-linear regime is full rank and only energy, mass and momentum preservation can be expected from low-rank techniques.
Again, we consider a setup similar to \cite{einkemmer2018low}, the domain is $\Omega = (0, 10 \pi) \times (-6, 6)$ and both directions are discretized with $n=128$ grid points and periodic boundary conditions. 
The initial value is
\begin{equation} \label{eq: two-stream instability}
f(0, x, v) = \frac{1}{2 \sqrt{2 \pi}} \left( e^{-(v-v_0)^2} + e^{-(v+v_0)^2} \right) \left( 1 + \alpha \cos(kx) \right),
\end{equation} 
where $\alpha = 10^{-3}$, $v_0=2.4$, and $k=1/5$.

\begin{figure}[ht!]
\centering
\includegraphics[width=0.8\textwidth]{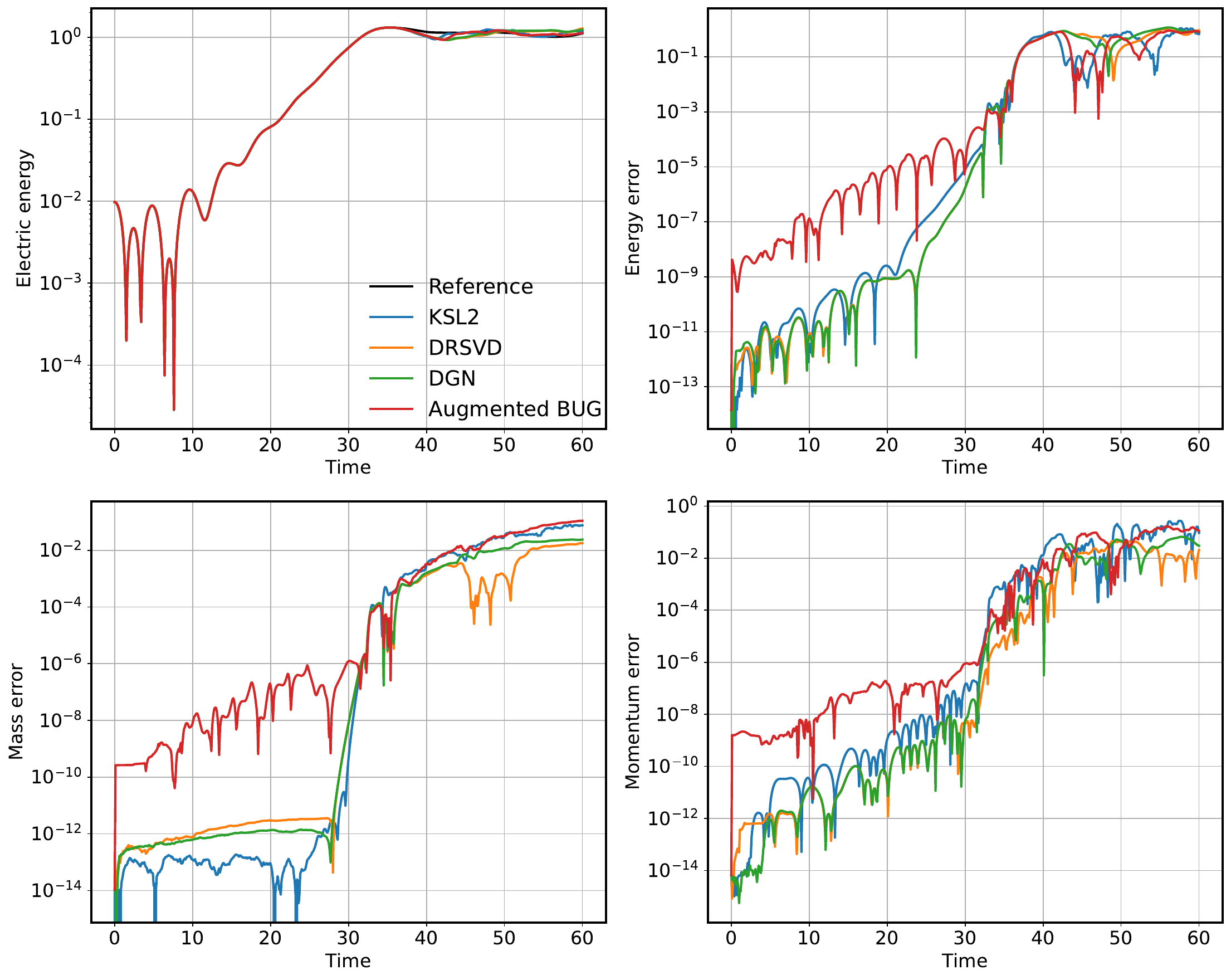}
\caption{Numerical simulations of the two-stream instability \eqref{eq: two-stream instability} with stepsize $h=0.1$. The DRSVD is performed with oversampling parameter $p=5$ and $q=1$ power iteration. The DGN is performed with oversampling parameters $p=5$, $\ell=0$ and $q=1$ power iteration.}
\label{fig: two-stream energy}
\end{figure}

\begin{figure}[H]
\centering
\includegraphics[width=\textwidth]{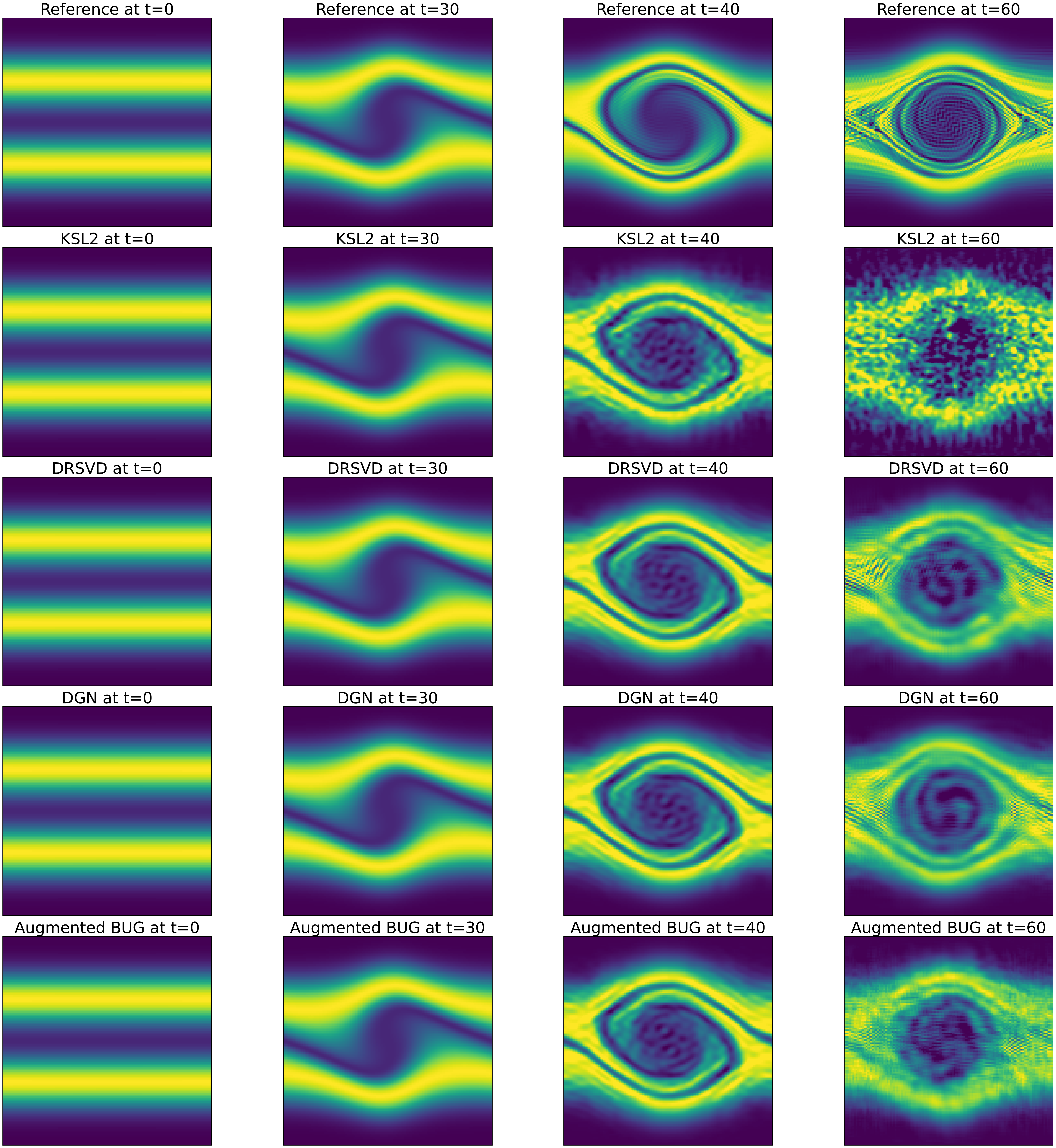}
\caption{Two-stream instability solution at several time steps ($t=0, 30, 40, 60$ from left to right). The solvers used are, from top to bottom: reference scipy RK45, projector-splitting, dynamical randomized SVD, dynamical generalized Nyström, augmented BUG. The color scale is shared for all plots.}
\label{fig: two-stream solutions}
\end{figure}

We compare the new methods to the projector-splitting and augmented BUG integrators in Figure \ref{fig: two-stream energy}.
In the linear regime, the new methods have errors comparable to the projector-splitting integrator, which is better than the augmented BUG integrator on this problem. In the non-linear regime, all methods keep the electric energy at a constant level.
In Figure~\ref{fig: two-stream solutions}, we see that the new methods seem to approximate better the geometry of the reference solution even though the non-linear regime does not have a low-rank structure. The solution obtained by the DGN algorithm looks particularly symmetrical compared to other low-rank methods.

\section*{Acknowledgements}

Thanks to Bart Vandereycken and Daniel Kressner for their enthusiasm about this project, which encouraged me to pursue the research and publish the paper as single author. I also want to express my gratitude to Bart Vandereycken for his valuable feedback on the early drafts. This work was supported by the SNSF under research project 192363.

\bibliographystyle{alpha}
\bibliography{references}

\newcommand{\etalchar}[1]{$^{#1}$}
\begin{thebibliography}{AAAAM{\etalchar{+}}21}

\bibitem[AAAAM{\etalchar{+}}21]{ahmadi2021randomized}
Salman Ahmadi-Asl, Stanislav Abukhovich, Maame~G Asante-Mensah, Andrzej
  Cichocki, Anh~Huy Phan, Tohishisa Tanaka, and Ivan Oseledets.
\newblock Randomized algorithms for computation of tucker decomposition and
  higher order svd (hosvd).
\newblock {\em IEEE Access}, 9:28684--28706, 2021.

\bibitem[AC72]{allen1972ground}
Samuel~Miller Allen and John~W Cahn.
\newblock Ground state structures in ordered binary alloys with second neighbor
  interactions.
\newblock {\em Acta Metallurgica}, 20(3):423--433, 1972.

\bibitem[AC73]{allen1973correction}
Samuel~M Allen and John~W Cahn.
\newblock A correction to the ground state of fcc binary ordered alloys with
  first and second neighbor pairwise interactions.
\newblock {\em Scripta Metallurgica}, 7(12):1261--1264, 1973.

\bibitem[AC25]{appelo2025robust}
Daniel Appel{\"o} and Yingda Cheng.
\newblock Robust implicit adaptive low rank time-stepping methods for matrix
  differential equations.
\newblock {\em Journal of Scientific Computing}, 102(3):81, 2025.

\bibitem[AL24]{ali2024dynamically}
Wael~H Ali and Pierre~FJ Lermusiaux.
\newblock Dynamically orthogonal narrow-angle parabolic equations for
  stochastic underwater sound propagation. part ii: Applications.
\newblock {\em The Journal of the Acoustical Society of America},
  155(1):656--672, 2024.

\bibitem[BCSK17]{babaee2017robust}
Hessam Babaee, Minseok Choi, Themistoklis~P Sapsis, and George~Em Karniadakis.
\newblock A robust bi-orthogonal/dynamically-orthogonal method using the
  covariance pseudo-inverse with application to stochastic flow problems.
\newblock {\em Journal of Computational Physics}, 344:303--319, 2017.

\bibitem[BEKK24]{baumann2024energy}
Lena Baumann, Lukas Einkemmer, Christian Klingenberg, and Jonas Kusch.
\newblock Energy stable and conservative dynamical low-rank approximation for
  the su--olson problem.
\newblock {\em SIAM Journal on Scientific Computing}, 46(2):B137--B158, 2024.

\bibitem[BFFN22]{billaud2022new}
Marie Billaud-Friess, Antonio Falc{\'o}, and Anthony Nouy.
\newblock A new splitting algorithm for dynamical low-rank approximation
  motivated by the fibre bundle structure of matrix manifolds.
\newblock {\em BIT Numerical Mathematics}, 62(2):387--408, 2022.

\bibitem[BFN17]{billaud2017dynamical}
Marie Billaud-Friess and Anthony Nouy.
\newblock Dynamical model reduction method for solving parameter-dependent
  dynamical systems.
\newblock {\em SIAM Journal on Scientific Computing}, 39(4):A1766--A1792, 2017.

\bibitem[Bit13]{bittencourt2013fundamentals}
Jos{\'e}~A Bittencourt.
\newblock {\em Fundamentals of plasma physics}.
\newblock Springer Science \& Business Media, 2013.

\bibitem[BYDW18]{batselier2018computing}
Kim Batselier, Wenjian Yu, Luca Daniel, and Ngai Wong.
\newblock Computing low-rank approximations of large-scale matrices with the
  tensor network randomized svd.
\newblock {\em SIAM Journal on Matrix Analysis and Applications},
  39(3):1221--1244, 2018.

\bibitem[CE22]{cassini2022efficient}
Fabio Cassini and Lukas Einkemmer.
\newblock Efficient 6d vlasov simulation using the dynamical low-rank framework
  ensign.
\newblock {\em Computer Physics Communications}, 280:108489, 2022.

\bibitem[CEKL24]{ceruti2024robust}
Gianluca Ceruti, Lukas Einkemmer, Jonas Kusch, and Christian Lubich.
\newblock A robust second-order low-rank bug integrator based on the midpoint
  rule.
\newblock {\em BIT Numerical Mathematics}, 64(3):30, 2024.

\bibitem[CGV23]{carrel2023low}
Benjamin Carrel, Martin~J Gander, and Bart Vandereycken.
\newblock Low-rank parareal: a low-rank parallel-in-time integrator.
\newblock {\em BIT Numerical Mathematics}, 63(1):13, 2023.

\bibitem[CH22]{coughlin2022efficient}
Jack Coughlin and Jingwei Hu.
\newblock Efficient dynamical low-rank approximation for the
  vlasov-amp{\`e}re-fokker-planck system.
\newblock {\em Journal of Computational Physics}, 470:111590, 2022.

\bibitem[CKL22]{ceruti2022rank}
Gianluca Ceruti, Jonas Kusch, and Christian Lubich.
\newblock A rank-adaptive robust integrator for dynamical low-rank
  approximation.
\newblock {\em BIT Numerical Mathematics}, 62(4):1149--1174, 2022.

\bibitem[CKL24]{ceruti2024parallel}
Gianluca Ceruti, Jonas Kusch, and Christian Lubich.
\newblock A parallel rank-adaptive integrator for dynamical low-rank
  approximation.
\newblock {\em SIAM Journal on Scientific Computing}, 46(3):B205--B228, 2024.

\bibitem[CKLS24]{ceruti2024parallelTTNs}
Gianluca Ceruti, Jonas Kusch, Christian Lubich, and Dominik Sulz.
\newblock A parallel basis update and galerkin integrator for tree tensor
  networks.
\newblock {\em arXiv preprint arXiv:2412.00858}, 2024.

\bibitem[CL20]{ceruti2020time}
Gianluca Ceruti and Christian Lubich.
\newblock Time integration of symmetric and anti-symmetric low-rank matrices
  and tucker tensors.
\newblock {\em BIT Numerical Mathematics}, 60:591--614, 2020.

\bibitem[CL22]{ceruti2022unconventional}
Gianluca Ceruti and Christian Lubich.
\newblock An unconventional robust integrator for dynamical low-rank
  approximation.
\newblock {\em BIT Numerical Mathematics}, 62(1):23--44, 2022.

\bibitem[CLS23]{ceruti2023rank}
Gianluca Ceruti, Christian Lubich, and Dominik Sulz.
\newblock Rank-adaptive time integration of tree tensor networks.
\newblock {\em SIAM Journal on Numerical Analysis}, 61(1):194--222, 2023.

\bibitem[CV23]{carrel2023projected}
Benjamin Carrel and Bart Vandereycken.
\newblock Projected exponential methods for stiff dynamical low-rank
  approximation problems.
\newblock {\em arXiv preprint arXiv:2312.00172}, 2023.

\bibitem[CW09]{clarkson2009numerical}
Kenneth~L Clarkson and David~P Woodruff.
\newblock Numerical linear algebra in the streaming model.
\newblock In {\em Proceedings of the forty-first annual ACM symposium on Theory
  of computing}, pages 205--214, 2009.

\bibitem[DPN{\etalchar{+}}23]{donello2023oblique}
M~Donello, G~Palkar, MH~Naderi, DC~Del Rey~Fern{\'a}ndez, and H~Babaee.
\newblock Oblique projection for scalable rank-adaptive reduced-order modelling
  of nonlinear stochastic partial differential equations with time-dependent
  bases.
\newblock {\em Proceedings of the Royal Society A}, 479(2278):20230320, 2023.

\bibitem[Ein19]{einkemmer2019low}
Lukas Einkemmer.
\newblock A low-rank algorithm for weakly compressible flow.
\newblock {\em SIAM Journal on Scientific Computing}, 41(5):A2795--A2814, 2019.

\bibitem[Ein24]{einkemmer2024accelerating}
Lukas Einkemmer.
\newblock Accelerating the simulation of kinetic shear alfv{\'e}n waves with a
  dynamical low-rank approximation.
\newblock {\em Journal of Computational Physics}, 501:112757, 2024.

\bibitem[EJ21]{einkemmer2021mass}
Lukas Einkemmer and Ilon Joseph.
\newblock A mass, momentum, and energy conservative dynamical low-rank scheme
  for the vlasov equation.
\newblock {\em Journal of Computational Physics}, 443:110495, 2021.

\bibitem[EKS23]{einkemmer2023conservation}
Lukas Einkemmer, Jonas Kusch, and Steffen Schotth{\"o}fer.
\newblock Conservation properties of the augmented basis update \& galerkin
  integrator for kinetic problems.
\newblock {\em arXiv preprint arXiv:2311.06399}, 2023.

\bibitem[EL18]{einkemmer2018low}
Lukas Einkemmer and Christian Lubich.
\newblock A low-rank projector-splitting integrator for the vlasov--poisson
  equation.
\newblock {\em SIAM Journal on Scientific Computing}, 40(5):B1330--B1360, 2018.

\bibitem[EMP24]{einkemmer2024hierarchical}
Lukas Einkemmer, Julian Mangott, and Martina Prugger.
\newblock A hierarchical dynamical low-rank algorithm for the stochastic
  description of large reaction networks.
\newblock {\em arXiv preprint arXiv:2407.11792}, 2024.

\bibitem[EOP20]{einkemmer2020low}
Lukas Einkemmer, Alexander Ostermann, and Chiara Piazzola.
\newblock A low-rank projector-splitting integrator for the vlasov--maxwell
  equations with divergence correction.
\newblock {\em Journal of Computational Physics}, 403:109063, 2020.

\bibitem[EOS23]{einkemmer2023robust}
Lukas Einkemmer, Alexander Ostermann, and Carmela Scalone.
\newblock A robust and conservative dynamical low-rank algorithm.
\newblock {\em Journal of Computational Physics}, 484:112060, 2023.

\bibitem[FBCM04]{fowlkes2004spectral}
Charless Fowlkes, Serge Belongie, Fan Chung, and Jitendra Malik.
\newblock Spectral grouping using the nystrom method.
\newblock {\em IEEE transactions on pattern analysis and machine intelligence},
  26(2):214--225, 2004.

\bibitem[FL18]{feppon2018dynamically}
Florian Feppon and Pierre~FJ Lermusiaux.
\newblock Dynamically orthogonal numerical schemes for efficient stochastic
  advection and lagrangian transport.
\newblock {\em Siam Review}, 60(3):595--625, 2018.

\bibitem[FLLN24]{feischl2024regularized}
Michael Feischl, Caroline Lasser, Christian Lubich, and J{\"o}rg Nick.
\newblock Regularized dynamical parametric approximation.
\newblock {\em arXiv preprint arXiv:2403.19234}, 2024.

\bibitem[FYL18]{feng2018faster}
Xu~Feng, Wenjian Yu, and Yaohang Li.
\newblock Faster matrix completion using randomized svd.
\newblock In {\em 2018 IEEE 30th International conference on tools with
  artificial intelligence (ICTAI)}, pages 608--615. IEEE, 2018.

\bibitem[GB24]{ghahremani2024cross}
Behzad Ghahremani and Hessam Babaee.
\newblock Cross interpolation for solving high-dimensional dynamical systems on
  low-rank tucker and tensor train manifolds.
\newblock {\em Computer Methods in Applied Mechanics and Engineering},
  432:117385, 2024.

\bibitem[GGH{\etalchar{+}}22]{greenacre2022principal}
Michael Greenacre, Patrick~JF Groenen, Trevor Hastie, Alfonso~Iodice d’Enza,
  Angelos Markos, and Elena Tuzhilina.
\newblock Principal component analysis.
\newblock {\em Nature Reviews Methods Primers}, 2(1):100, 2022.

\bibitem[Git11]{gittens2011spectral}
Alex Gittens.
\newblock The spectral norm error of the naive nystrom extension.
\newblock {\em arXiv preprint arXiv:1110.5305}, 2011.

\bibitem[GM13]{gittens2013revisiting}
Alex Gittens and Michael Mahoney.
\newblock Revisiting the nystrom method for improved large-scale machine
  learning.
\newblock In {\em International Conference on Machine Learning}, pages
  567--575. PMLR, 2013.

\bibitem[GQ22]{guo2022low}
Wei Guo and Jing-Mei Qiu.
\newblock A low rank tensor representation of linear transport and nonlinear
  vlasov solutions and their associated flow maps.
\newblock {\em Journal of Computational Physics}, 458:111089, 2022.

\bibitem[GQ24]{guo2024conservative}
Wei Guo and Jing-Mei Qiu.
\newblock A conservative low rank tensor method for the vlasov dynamics.
\newblock {\em SIAM Journal on Scientific Computing}, 46(1):A232--A263, 2024.

\bibitem[HMT11]{halko2011finding}
Nathan Halko, Per-Gunnar Martinsson, and Joel~A Tropp.
\newblock Finding structure with randomness: Probabilistic algorithms for
  constructing approximate matrix decompositions.
\newblock {\em SIAM review}, 53(2):217--288, 2011.

\bibitem[JL14]{ji2014gpu}
Hao Ji and Yaohang Li.
\newblock Gpu accelerated randomized singular value decomposition and its
  application in image compression.
\newblock {\em Proc. of MSVESCC}, pages 39--45, 2014.

\bibitem[KEC23]{kusch2023stability}
Jonas Kusch, Lukas Einkemmer, and Gianluca Ceruti.
\newblock On the stability of robust dynamical low-rank approximations for
  hyperbolic problems.
\newblock {\em SIAM Journal on Scientific Computing}, 45(1):A1--A24, 2023.

\bibitem[KL07]{koch2007dynamical}
Othmar Koch and Christian Lubich.
\newblock Dynamical low-rank approximation.
\newblock {\em SIAM Journal on Matrix Analysis and Applications},
  29(2):434--454, 2007.

\bibitem[KL10]{koch2010dynamical}
Othmar Koch and Christian Lubich.
\newblock Dynamical tensor approximation.
\newblock {\em SIAM Journal on Matrix Analysis and Applications},
  31(5):2360--2375, 2010.

\bibitem[KLW16]{kieri2016discretized}
Emil Kieri, Christian Lubich, and Hanna Walach.
\newblock Discretized dynamical low-rank approximation in the presence of small
  singular values.
\newblock {\em SIAM Journal on Numerical Analysis}, 54(2):1020--1038, 2016.

\bibitem[KNZ24]{kazashi2024dynamical}
Yoshihito Kazashi, Fabio Nobile, and Fabio Zoccolan.
\newblock Dynamical low-rank approximation for stochastic differential
  equations.
\newblock {\em Mathematics of Computation}, 2024.

\bibitem[KS23]{kusch2023robust}
Jonas Kusch and Pia Stammer.
\newblock A robust collision source method for rank adaptive dynamical low-rank
  approximation in radiation therapy.
\newblock {\em ESAIM: Mathematical Modelling and Numerical Analysis},
  57(2):865--891, 2023.

\bibitem[KSW25]{kusch2025augmented}
Jonas Kusch, Steffen Schotth{\"o}fer, and Alexandra Walter.
\newblock An augmented backward-corrected projector splitting integrator for
  dynamical low-rank training.
\newblock {\em arXiv preprint arXiv:2502.03006}, 2025.

\bibitem[Kus24]{kusch2024second}
Jonas Kusch.
\newblock Second-order robust parallel integrators for dynamical low-rank
  approximation.
\newblock {\em arXiv preprint arXiv:2403.02834}, 2024.

\bibitem[KV19]{kieri2019projection}
Emil Kieri and Bart Vandereycken.
\newblock Projection methods for dynamical low-rank approximation of
  high-dimensional problems.
\newblock {\em Computational Methods in Applied Mathematics}, 19(1):73--92,
  2019.

\bibitem[LBKL14]{li2014large}
Mu~Li, Wei Bi, James~T Kwok, and Bao-Liang Lu.
\newblock Large-scale nystr{\"o}m kernel matrix approximation using randomized
  svd.
\newblock {\em IEEE transactions on neural networks and learning systems},
  26(1):152--164, 2014.

\bibitem[LCK24]{lam2024randomized}
Hei~Yin Lam, Gianluca Ceruti, and Daniel Kressner.
\newblock Randomized low-rank runge-kutta methods.
\newblock {\em arXiv preprint arXiv:2409.06384}, 2024.

\bibitem[LL03]{lee2003smooth}
John~M Lee and John~M Lee.
\newblock {\em Smooth manifolds}.
\newblock Springer, 2003.

\bibitem[LO14]{lubich2014projector}
Christian Lubich and Ivan~V Oseledets.
\newblock A projector-splitting integrator for dynamical low-rank
  approximation.
\newblock {\em BIT Numerical Mathematics}, 54(1):171--188, 2014.

\bibitem[LRSV13]{lubich2013dynamical}
Christian Lubich, Thorsten Rohwedder, Reinhold Schneider, and Bart
  Vandereycken.
\newblock Dynamical approximation by hierarchical tucker and tensor-train
  tensors.
\newblock {\em SIAM Journal on Matrix Analysis and Applications},
  34(2):470--494, 2013.

\bibitem[LVW18]{lubich2018time}
Christian Lubich, Bart Vandereycken, and Hanna Walach.
\newblock Time integration of rank-constrained tucker tensors.
\newblock {\em SIAM Journal on Numerical Analysis}, 56(3):1273--1290, 2018.

\bibitem[Mar19]{martinsson2019randomized}
Per-Gunnar Martinsson.
\newblock Randomized methods for matrix computations.
\newblock {\em The Mathematics of Data}, 25(4):187--231, 2019.

\bibitem[MT20]{martinsson2020randomized}
Per-Gunnar Martinsson and Joel~A Tropp.
\newblock Randomized numerical linear algebra: Foundations and algorithms.
\newblock {\em Acta Numerica}, 29:403--572, 2020.

\bibitem[NAB24]{naderi2024cur}
Mohammad~Hossein Naderi, Sara Akhavan, and Hessam Babaee.
\newblock Cur for implicit time integration of random partial differential
  equations on low-rank matrix manifolds.
\newblock {\em arXiv preprint arXiv:2408.16591}, 2024.

\bibitem[Nak20]{nakatsukasa2020fast}
Yuji Nakatsukasa.
\newblock Fast and stable randomized low-rank matrix approximation.
\newblock {\em arXiv preprint arXiv:2009.11392}, 2020.

\bibitem[NL08]{nonnenmacher2008dynamical}
Achim Nonnenmacher and Christian Lubich.
\newblock Dynamical low-rank approximation: applications and numerical
  experiments.
\newblock {\em Mathematics and Computers in Simulation}, 79(4):1346--1357,
  2008.

\bibitem[NQE25]{nakao2025reduced}
Joseph Nakao, Jing-Mei Qiu, and Lukas Einkemmer.
\newblock Reduced augmentation implicit low-rank (rail) integrators for
  advection-diffusion and fokker--planck models.
\newblock {\em SIAM Journal on Scientific Computing}, 47(2):A1145--A1169, 2025.

\bibitem[NR25]{nobile2025robust}
Fabio Nobile and S{\'e}bastien Riffaud.
\newblock Robust high-order low-rank bug integrators based on explicit
  runge-kutta methods.
\newblock {\em arXiv preprint arXiv:2502.07040}, 2025.

\bibitem[NT23]{nobile2023error}
Fabio Nobile and Thomas~Trigo Trindade.
\newblock Error estimates for supg-stabilised dynamical low rank
  approximations.
\newblock In {\em European Conference on Numerical Mathematics and Advanced
  Applications}, pages 438--447. Springer, 2023.

\bibitem[NT25]{nobile2025petrov}
Fabio Nobile and Thomas~Trigo Trindade.
\newblock Petrov--galerkin dynamical low rank approximation: Supg stabilisation
  of advection-dominated problems.
\newblock {\em Computer Methods in Applied Mechanics and Engineering},
  433:117495, 2025.

\bibitem[Nys30]{nystrom1930praktische}
Evert~J Nystr{\"o}m.
\newblock {\"U}ber die praktische aufl{\"o}sung von integralgleichungen mit
  anwendungen auf randwertaufgaben.
\newblock 1930.

\bibitem[OPW19]{ostermann2019convergence}
Alexander Ostermann, Chiara Piazzola, and Hanna Walach.
\newblock Convergence of a low-rank lie--trotter splitting for stiff matrix
  differential equations.
\newblock {\em SIAM Journal on Numerical Analysis}, 57(4):1947--1966, 2019.

\bibitem[PN25]{park2025low}
Taejun Park and Yuji Nakatsukasa.
\newblock Low-rank approximation of parameter-dependent matrices via cur
  decomposition.
\newblock {\em SIAM Journal on Scientific Computing}, 47(3):A1858--A1887, 2025.

\bibitem[RST10]{rokhlin2010randomized}
Vladimir Rokhlin, Arthur Szlam, and Mark Tygert.
\newblock A randomized algorithm for principal component analysis.
\newblock {\em SIAM Journal on Matrix Analysis and Applications},
  31(3):1100--1124, 2010.

\bibitem[RV23]{rodgers2023implicit}
Abram Rodgers and Daniele Venturi.
\newblock Implicit integration of nonlinear evolution equations on tensor
  manifolds.
\newblock {\em Journal of Scientific Computing}, 97(2):33, 2023.

\bibitem[Sch22]{schmid2022dynamic}
Peter~J Schmid.
\newblock Dynamic mode decomposition and its variants.
\newblock {\em Annual Review of Fluid Mechanics}, 54(1):225--254, 2022.

\bibitem[SEKM25]{scalone2025multi}
Carmela Scalone, Lukas Einkemmer, Jonas Kusch, and Ryan~G McClarren.
\newblock A multi-fidelity adaptive dynamical low-rank based optimization
  algorithm for fission criticality problems.
\newblock {\em Journal of Scientific Computing}, 104(1):1--18, 2025.

\bibitem[SHNT23]{schmidt2023rank}
Jonathan Schmidt, Philipp Hennig, J{\"o}rg Nick, and Filip Tronarp.
\newblock The rank-reduced kalman filter: Approximate dynamical-low-rank
  filtering in high dimensions.
\newblock {\em Advances in Neural Information Processing Systems},
  36:61364--61376, 2023.

\bibitem[SHvBW21]{saibaba2021randomized}
Arvind~K Saibaba, Joseph Hart, and Bart van Bloemen~Waanders.
\newblock Randomized algorithms for generalized singular value decomposition
  with application to sensitivity analysis.
\newblock {\em Numerical linear algebra with applications}, 28(4):e2364, 2021.

\bibitem[SL24]{schotthofer2024federated}
Steffen Schotth{\"o}fer and M~Paul Laiu.
\newblock Federated dynamical low-rank training with global loss convergence
  guarantees.
\newblock {\em arXiv preprint arXiv:2406.17887}, 2024.

\bibitem[SLC{\etalchar{+}}24]{sulz2024numerical}
Dominik Sulz, Christian Lubich, Gianluca Ceruti, Igor Lesanovsky, and Federico
  Carollo.
\newblock Numerical simulation of long-range open quantum many-body dynamics
  with tree tensor networks.
\newblock {\em Physical Review A}, 109(2):022420, 2024.

\bibitem[SS02]{scholkopf2002learning}
Bernhard Sch{\"o}lkopf and Alexander~J Smola.
\newblock {\em Learning with kernels: support vector machines, regularization,
  optimization, and beyond}.
\newblock MIT press, 2002.

\bibitem[SV24]{sutti2024implicit}
Marco Sutti and Bart Vandereycken.
\newblock Implicit low-rank riemannian schemes for the time integration of
  stiff partial differential equations.
\newblock {\em Journal of Scientific Computing}, 101(1):3, 2024.

\bibitem[SYS25]{schotthofer2025dynamical}
Steffen Schotth{\"o}fer, H~Lexie Yang, and Stefan Schnake.
\newblock Dynamical low-rank compression of neural networks with robustness
  under adversarial attacks.
\newblock {\em arXiv preprint arXiv:2505.08022}, 2025.

\bibitem[SZC{\etalchar{+}}24]{schotthofer2024geolora}
Steffen Schotth{\"o}fer, Emanuele Zangrando, Gianluca Ceruti, Francesco
  Tudisco, and Jonas Kusch.
\newblock Geolora: Geometric integration for parameter efficient fine-tuning.
\newblock {\em arXiv preprint arXiv:2410.18720}, 2024.

\bibitem[SZCT23]{savostianova2023robust}
Dayana Savostianova, Emanuele Zangrando, Gianluca Ceruti, and Francesco
  Tudisco.
\newblock Robust low-rank training via approximate orthonormal constraints.
\newblock {\em Advances in Neural Information Processing Systems},
  36:66064--66083, 2023.

\bibitem[SZK{\etalchar{+}}22]{schotthofer2022low}
Steffen Schotth{\"o}fer, Emanuele Zangrando, Jonas Kusch, Gianluca Ceruti, and
  Francesco Tudisco.
\newblock Low-rank lottery tickets: finding efficient low-rank neural networks
  via matrix differential equations.
\newblock {\em Advances in Neural Information Processing Systems},
  35:20051--20063, 2022.

\bibitem[TYUC17]{tropp2017practical}
Joel~A Tropp, Alp Yurtsever, Madeleine Udell, and Volkan Cevher.
\newblock Practical sketching algorithms for low-rank matrix approximation.
\newblock {\em SIAM Journal on Matrix Analysis and Applications},
  38(4):1454--1485, 2017.

\bibitem[UV20]{uschmajew2020geometric}
Andr{\'e} Uschmajew and Bart Vandereycken.
\newblock {\em Geometric methods on low-rank matrix and tensor manifolds}.
\newblock Springer, 2020.

\bibitem[UZ24]{uschmajew2024dynamical}
Andr{\'e} Uschmajew and Andreas Zeiser.
\newblock Dynamical low-rank approximation of the vlasov--poisson equation with
  piecewise linear spatial boundary.
\newblock {\em BIT Numerical Mathematics}, 64(2):19, 2024.

\bibitem[W{\etalchar{+}}14]{woodruff2014sketching}
David~P Woodruff et~al.
\newblock Sketching as a tool for numerical linear algebra.
\newblock {\em Foundations and Trends{\textregistered} in Theoretical Computer
  Science}, 10(1--2):1--157, 2014.

\bibitem[WLC23]{wen2023accelerated}
You-Wei Wen, Kexin Li, and Hefeng Chen.
\newblock Accelerated matrix completion algorithm using continuation strategy
  and randomized svd.
\newblock {\em Journal of Computational and Applied Mathematics}, 429:115215,
  2023.

\bibitem[WLRT08]{woolfe2008fast}
Franco Woolfe, Edo Liberty, Vladimir Rokhlin, and Mark Tygert.
\newblock A fast randomized algorithm for the approximation of matrices.
\newblock {\em Applied and Computational Harmonic Analysis}, 25(3):335--366,
  2008.

\bibitem[WS00]{williams2000using}
Christopher Williams and Matthias Seeger.
\newblock Using the nystr{\"o}m method to speed up kernel machines.
\newblock {\em Advances in neural information processing systems}, 13, 2000.

\bibitem[XZ13]{xiang2013regularization}
Hua Xiang and Jun Zou.
\newblock Regularization with randomized svd for large-scale discrete inverse
  problems.
\newblock {\em Inverse Problems}, 29(8):085008, 2013.

\bibitem[YEHS23]{yin2023semi}
Peimeng Yin, Eirik Endeve, Cory~D Hauck, and Stefan~R Schnake.
\newblock A semi-implicit dynamical low-rank discontinuous galerkin method for
  space homogeneous kinetic equations. part i: emission and absorption.
\newblock {\em arXiv preprint arXiv:2308.05914}, 2023.

\end{thebibliography}

\end{document}